\documentclass[rewiew]{article}

\usepackage{amsfonts, amsmath, amsthm, amssymb}

\usepackage{mathtools}
\usepackage{makecell}
\usepackage{pgf,tikz}
\usetikzlibrary{arrows}

\usepackage{dsfont}
\usepackage{secdot}

\usepackage{algorithm}
\usepackage{algorithmic}
\usepackage{float}
\usepackage{comment}

\usepackage{multicol}
\usepackage{lipsum}
\usepackage{mwe}

\usepackage{subcaption}

\newtheorem{lemma}{Lemma}

\newtheorem{definition}{Definition}

\makeatletter
\renewcommand{\ALG@name}{Algorithm}

\allowdisplaybreaks

\graphicspath{{pictures/}} 

\tikzset{every picture/.style={line width=0.7pt}}
\tikzset{every node/.style={scale=1.5}}

\title{Coverings of planar and three-dimensional sets with subsets of smaller diameter \footnote{The work was supported by the program “Leading Scientific Schools” under grant NSh-775.2022.1.1.}}

\author{A.\,D.\,Tolmachev\footnote{Graduate student, Moscow Institute of Physics and Technology}, D.\,S.\,Protasov\footnote{Graduate student, Moscow Institute of Physics and Technology}, V.\,A.\,Voronov\footnote{PhD, Researcher, Caucasus Mathematical Center of Adyghe State University, Maikop; Moscow Institute of Physics and Technology} \\{\small\textit{v-vor@yandex.ru}}}

\begin{document}

\maketitle

\begin{abstract}
Quantitative estimates related to the classical Borsuk problem of splitting set in Euclidean space into subsets of smaller diameter are considered. For a given $k$ there is a minimal diameter of subsets at which there exists a covering with $k$ subsets of any planar set of unit diameter.  In order to find an upper estimate of the minimal diameter we propose an algorithm for finding sub-optimal partitions. In the cases  $10 \leqslant k \leqslant 17$ some upper and lower estimates of the minimal diameter are improved. Another result is that any set $M \subset \mathbb{R}^3$ of a unit diameter can be partitioned into four subsets of a diameter not greater than $0.966$.
\end{abstract}




\section{Introduction}

 This paper presents some results related to the classical Borsuk problem  on partitioning of sets in $\mathbb{R}^n$ into parts of smaller diameter \cite{borsuk,borsuk1,borsuk2,borsuk3}, and also to the Nelson--Erdős--Hadwiger problem  on the chromatic number of Euclidean space \cite{chrom,chrom1,chrom2,chrom3,chrom4,chrom5,chrom6}. Most of the article is devoted to the planar case.

Let $F$ be a bounded set in the plane, and $k \in \mathbb{N}$. We denote by $d_k(F)$ the greatest lower bound of the set of positive real numbers $x$ with the property that $F$ can be covered by $k$ sets $F_1, F_2, \ldots , F_k$ whose diameters are at most $x$, that is,
 $$d_k(F) = \inf\{x \in \mathbb{R}^{+} : \exists F_{1}, \ldots, F_{k}: F \subseteq F_{1} \cup \ldots \cup F_{k}, \; \forall i \:  \operatorname{diam}(F_{i}) \leqslant x \}.$$

 In other words, we want to choose optimal coverings consisting of the smallest possible diameter sets among all possible coverings of the set $F$. In addition, the value $d_k(F)$ does not change, if, without loss of generality, we require the sets to be convex and closed. Indeed, it is easy to see that the diameter of the closure of the convex hull for any set $F_i$ from the covering coincides with the diameter of $F_i$.
 
For every positive number $k$ we consider the values $d_k = \sup \; d_k(F)$, where the suprema are taken over all sets $F$ of unit diameter on the plane. It follows from the remark above that the sequence $d_k$ is nonincreasing.

Motivated by the classical Borsuk problem  \cite{borsuk,borsuk1,borsuk2,borsuk3}, many specialists have evaluated elements of this sequence. Over the years, H. Lenz (see \cite{Lenz}), M. Dembinski and M. Lassak (see \cite{Lassak}), V. Filimonov (see \cite{Filimonov}), D. Belov and N. Aleksandrov (see \cite{Belov}), V. Koval (see \cite{Koval}) estimated values $d_k$ for various values of $k$. A theoretically feasible, but extremely time-consuming approach to this problem and its generalizations was proposed in \cite{zong2021borsuk}. Moreover, Yanlu Lian and Senlin Wu explored such values for some Banach spaces (see \cite{lian2021divide}). In our previous work (see \cite{Doklady}) we significantly improved some of previous upper bounds of the quantities $d_k$. In this paper we prove new lower and upper bounds for the elements in the sequence $d_k$. Moreover, we proposed  new approach to improving upper bounds of the values $d_k$.

We present our new results in the following sections, but here we suggest additional important definitions for these theorems and their proofs. In this paper, using techniques to construct universal covering sets and systems, we prove some upper bounds for the elements of the sequence $d_k$.

Note that both infinitesimal local improvements to these partitions are possible, as well as improvements based on the extension of the covering system. Of course, this approach does not allow us to obtain exact values of $d_{k}$.

\vspace{\baselineskip}
 
 \begin{definition}
  A set $\Omega \subset \mathbb R^2$ is called \textit{universal covering set} if every planar set $F$ of unit diameter can be completely covered by $\Omega$ (that is,
there exists a planar set $\Omega'$ congruent to $\Omega$ such that $F \subset \Omega'$ ).
\end{definition}
 
 In 1920, in \cite{Pal} J. Pal proved that a regular hexagon with edge length $\frac{1}{\sqrt{3}}$ is a universal covering set. We denote this regular hexagon by $\Omega$. Next, we define a universal covering system.
 
 \begin{definition}
 System of sets $S = \{ \Omega_{\alpha} \}_{\alpha \in I}$ is called a \textit{universal covering system} if every planar set  $F$ of unit diameter can be completely covered by one of sets $\Omega_{\alpha}$. Here $I$ is a (possibly infinite) set of indices.
 \end{definition}

\section{Main results}

 In this part of the paper we show the table of improved results for the first 17 elements of the sequence $d_k$. In Table 1, the column titled as ``comment'' indicates by how many percent the gap between the upper and lower bounds has decreased as a result of the improvements proposed in this paper. The column titled as ``UCS'' presents the universal covering system used to prove the indicated upper bound on $d_k$ for the corresponding value of $k$. Let us denote by $\underline{d}_k^{old}$, $\underline{d}_k^{new}$,  $\overline{d}_k^{old}$ and $\overline{d}_k^{new}$, respectively, the previously known and the new value of the lower bound, and similarly for the upper bound.
 
 All constants in this table are specified with four decimal places.

 \setlength{\tabcolsep}{8pt}
\renewcommand{\arraystretch}{0.95}

\begin{center}
\begin{table}[H]
\caption{Old and new estimates}
\label{results1}
 \begin{tabular}{ |c|c|c||c|c||c|c|  }
 \hline

 $k$ & $\underline{d}_k^{new}$ &
 $\underline{d}_k^{old}$ &  $\overline{d}_k^{old}$ & $\overline{d}_k^{new}$ & Comment & UCS \\
 \hline\hline
  1 & - & 1.0000 & 1.0000 & - & tight & -\\ 
  \hline
  2 & - & 1.0000 & 1.0000 & - & tight & -\\
  \hline
  3 & - & 0.8660 & 0.8660 & - & tight & -\\
  \hline
  4 & - & 0.7071 & 0.7071 & - & tight & -\\
  \hline
  5 & - & 0.5877 & 0.5953 & - & - &  - \\
  \hline
  6 & - & 0.5051 &  0.5343 & - & - & -\\
  \hline
  7 & - & 0.5000 & 0.5000 & - & tight & -\\
  \hline
  8 & - & 0.4338 & 0.4456 & - &- & -\\
  \hline
  9 & - & 0.3826 & 0.4047 & - & - & -\\
  \hline
  10 & 0.3665 & 0.3420 & 0.4012 & 0.3896 & 61\%  & $S_{10}$ \\
  \hline
  11 & 0.3535 & 0.3333 & 0.3942 & 0.3732 & 68\% & $S_{10}$ \\
  \hline

  12 & 0.3420  & 0.3333 & 0.3660 & 0.3532 & 66\% & $S_{10}$ \\
  \hline
  13 & - & 0.3333 & 0.3550 & 0.3419 & 60\% & $S_{10}$ \\
  \hline
  14 & - & 0.3090 & 0.3324 & 0.3263 & 26\% & $S_{10}$ \\
  \hline
  15 & - & 0.2928 & 0.3226 & 0.3130 & 32\% & $S_{10}$ \\
  \hline
  16 & - & 0.2817 & 0.3191 & 0.3035 & 42\% & $S_{10}$ \\
  \hline
  17 & - & 0.2701 & 0.3010 & 0.2967 & 14\% & $S_{10}$ \\
  \hline
 
\end{tabular}
\end{table}
\end{center}

\section{Upper bounds}

\subsection {Improvements to the upper estimates and the construction of the UCS}

To prove that $d_k \leqslant \rho$, where $\rho$ is some fixed number, we consider some universal covering system $S$ and divide each of the covering sets into $k$ parts. The diameter of each part does not exceed $\rho$.

From the introduction, we know that any set of diameter 1 can be covered by a $ \Omega$, i.e regular hexagon with a unit distance between the opposite sides.

\begin{figure}[htb!]
\center{\includegraphics[scale=0.4]{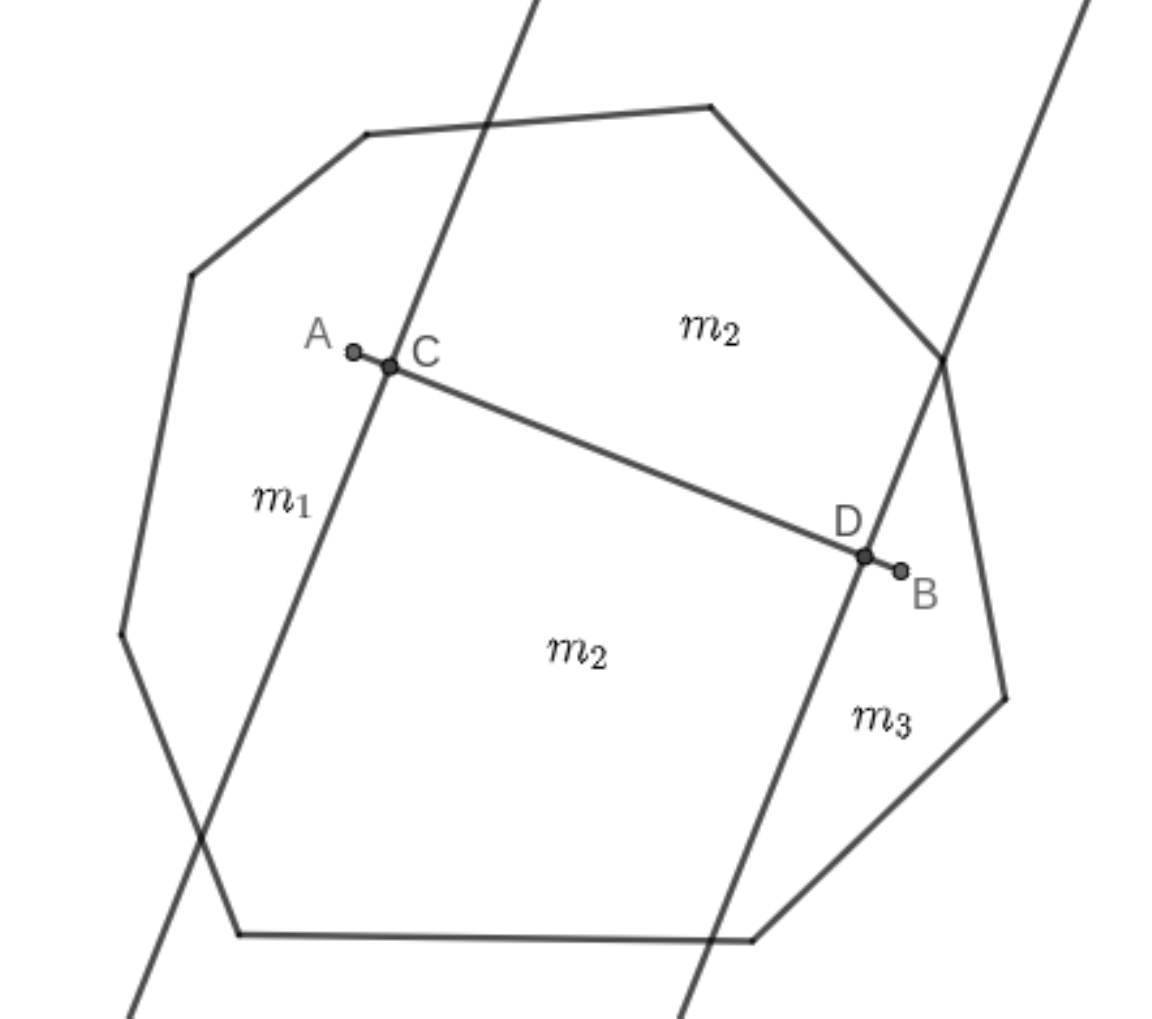}}
\caption{The set $M$ and its division into the sets $m_1$, $m_2$, $m_3$}
\label{cut1}
\end{figure}

\begin{lemma}
 Let $\sigma$ be a covering system. Suppose that there exists a set $M \in \sigma$ and points $A, B \in M$ such that the length of segment $AB$ is greater than 1. Denote by $\delta = \frac{|AB| - 1}{2}$. Perpendicular lines drawn to the segment $AB$ at a distance of $\frac{\delta}{2}$ from its ends divide the set $M$ into three sets $m_1$, $m_2$, $m_3$ (in Figure 1, the perpendiculars are drawn at points $C$ and $D$). We will refer the perpendiculars themselves only to the set $m_2$, so sets $m_1$ and $m_3$ will not contain them. Denote by $M_1 = m_1 \cup m_2$, $M_2 = m_2 \cup m_3$. 

Then $ \sigma \backslash \{M\}  \cup \{M_1, M_2\}$ will also be a covering system. 
\end{lemma}

\textbf {Proof:} We want to show that if a set $F$ can be covered by the set $M$, then $F$ can be covered by the set $M_1$ or the set $M_2$. Let's assume that when a certain set $N$ is covered by the set $M$, there is at least one point of the set $N$ in $m_1$. Then $m_3$ cannot have points from $N$, because the distance between any point from $m_1$ and a point from $m_3$ is strictly greater than 1. So we get a contradiction with the diameter of the set $N$ is not more than 1. This contradiction proves that the set $M_1$ will cover the set $N$. If $m_1$ does not contain points from $N$, then $M_2$ will cover $N$. This statement proves the Lemma 1. \qed

Now, using the lemma, we pass from one set $\Omega$ to a system of two covering sets. Each of the main diagonals of $ \Omega$ has length $\frac{2}{\sqrt{3}}$, which is greater than one, so according to the lemma from the hexagon, you can cut off the corresponding parts (in this case, triangles) on each of the three diagonals of the hexagon.

\begin{figure}[htb!]
\center{\includegraphics[scale=0.4]{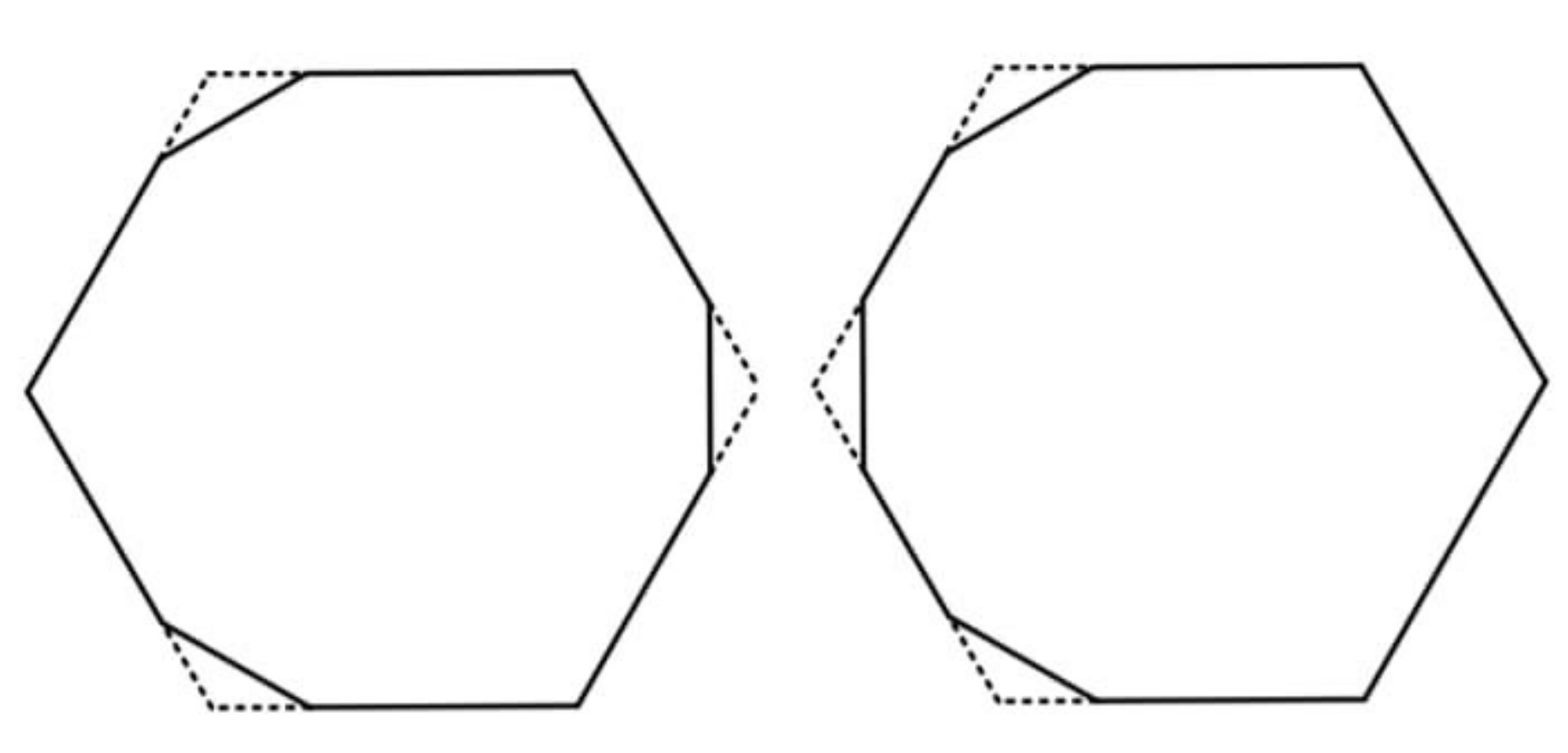}}
\caption{$\Omega_1, \Omega_2$ (the dotted line marks the parts cut off from $\Omega$)}
\label{cut_hex}
\end{figure}

After eliminating the congruent ones, we get a universal covering system containing two sets. In the first case, the vertices from which the triangles are cut off go through one (let's call it $\Omega_1$). In the second case, the vertices from which the triangles are cut off go in a row (let's call it $\Omega_2$).

\begin{figure}[!htb]
\center{\includegraphics[scale=0.25]{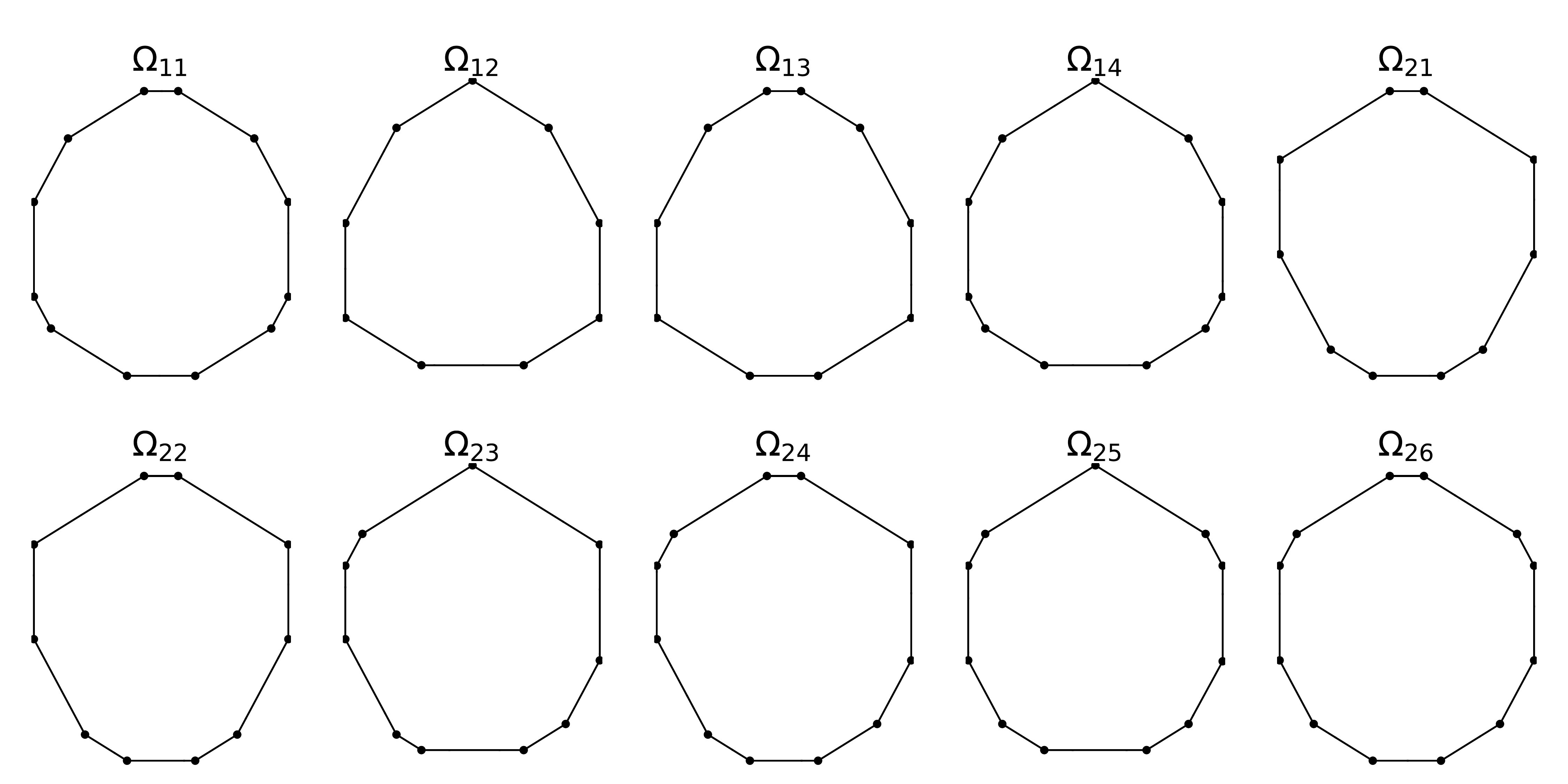}}
\caption{Shapes used in UCS $S_{10}$}
\label{UCS10}
\end{figure}

Note that after cutting off the triangles, the remaining segment of the diagonal (equal to $\frac{1}{2} + \frac{1}{\sqrt{3}}$) is still greater than $1$. Therefore, the lemma allows us to cut off from $\Omega_1, \Omega_2$ other parts, reducing the length of the diagonals by half the remaining distance (by $h = \frac{1}{4} + \frac{1}{2\sqrt{3}}$). There are two clipping options for each diagonal. In the first case, we cut off the triangle of height $h$, and in the other we cut off the trapezoid of height $h$. Excluding the congruent ones, we get a universal covering system consisting of 10 sets (4 sets are obtained in the results of cutting off pieces from $\Omega_1$, the other 6 in the results of cutting off pieces from $\Omega_2$). 
We denote by $S_{10} = \{\Omega_{11}, \Omega_{12}, \Omega_{13}, \Omega_{14}, \Omega_{21}, \Omega_{22}, \Omega_{23}, \Omega_{24}, \Omega_{25}, \Omega_{26}\}$ the constructed universal covering system of 10 sets. Using this UCS, new upper bounds are obtained for $d_{10}, \dots, d_{17}$.

All the necessary partitions of sets from the UCS are given in \cite{github} and obtained using an optimization algorithm.

The proposed method can be applied without any significant changes to the search for partitions in higher dimensions. Let us denote by $d_{3,4}$ the minimal diameter of a part for which there exists a partition of any three-dimensional set of unit diameter into four parts of a given diameter. The previously known estimate $d_{3,4} \leqslant 0.98$ was obtained by L. Evdokimov by partitioning a truncated rhombic dodecahedron $\Omega_{TRD}$ (Fig. \ref{fig:TRD1}). V.~V.~Makeev proved that this polyhedron is a universal cover in the three-dimensional case \cite{makeev2000}.  Using the optimization algorithm described in the next section we obtained a slightly better estimate  $$d_{3,4} \leqslant 0.966.$$ Note that here we consider a \emph{covering} of $\Omega_{TRD}$ by several convex polyhedrons, which cannot be directly reduced to \emph{partition} into convex polyhedrons. 

If we consider only the partitions of $\Omega_{TRD}$ into four convex polyhedrons of smaller diameter by six planes passing through some common point, the estimate is slightly worse, $d_{3,4} \leqslant 0.9755$, close to the result presented in \cite{makeev2000}. These partitions are shown in Fig. \ref{fig:TRD1}, Fig. \ref{fig:TRD2}. The coordinates of the vertices are available in the repository \cite{github}.

\subsection{Description of the algorithm}

The idea of the proposed algorithm lies in multiple generation of some initial partition into polygons and subsequent minimization of the maximum diameter of obtained parts.  We assume that when solving the optimization problem, we can consider the structure of the partition to be unchanged. Thus we are talking about finding the local minimum of some piecewise smooth nonconvex function under linear constraints.

The Adam algorithm, which was proposed in  \cite{adam} and is now widely used in machine learning problems, is used to find a local minimum.  The Adam algorithm is one of the extensions of stochastic gradient descent. Its high convergence rate and stability are achieved by adaptive learning rate selection for each parameter based on the mean and the variance of the gradient. Most theoretical results for algorithms of this type are obtained under the assumptions of differentiability and convexity of the minimizing function, but in machine learning problems the objective function is usually non-smooth, for example, when training a neural network with ReLU activation function \cite{adam}.

Note that in our case the minimized function is not stochastic. We use only those properties of the Adam algorithm that allow us to efficiently solve non-smooth high-dimensional optimization problems. In this paper, we do not prove any convergence statements and present only numerically found local minima.

In numerical calculations, we used a penalty method and an implementation of the Adam algorithm in the PyTorch package \cite{adam}.

Let some initial partition of an $m$-gon $\Omega$ into polygons $\Omega = F_1 \cup F_2 \cup \dots \cup F_k$ be chosen and $X = \{x_1, \dots , x_r \}$ be the vertex set of the partition.  The vertices of the polygon $F_i$ are points $x_j$, $j \in \mathcal{J}_i$. The condition $(x,c)+b=0$, $\|c\|=1$ specifies that the point $x$ belongs to the line given by the normal vector $c$  and the coefficient $b$. For an interior point $y \in \operatorname{Int} \Omega$ the value of the linear function $f(y) = (y,c)+b >0$ is the distance to the line. We assume that the belonging of points to lines bounding $\Omega$ is given for some set of pairs of indices $E \subseteq \{1, \dots r \} \times \{1, \dots , m\}$. Then the search for a local minimum in the problem 
\begin{align}
\label{opt1}
\varphi(X) = \max_{i} \max_{p,q \in \mathcal{J}_i} \| x_p - x_q \| \rightarrow \min, \\
(x_s, c_t) +b_t = 0, \quad (s,t) \in E \nonumber
\end{align}

\noindent can be performed by standard methods for non-smooth optimization problems. Generally speaking, we should also require the conditions
\[
   (x_s, c_t) +b_t \geq 0, \quad 1 \leq s \leq r, \quad 1 \leq t \leq m,
\]
meaning that each vertex $x_s$ belongs to $\Omega$. But if we assume that $x_s$ change rather little during partition optimization, we may not introduce these conditions, and check them after the local minimum has been found.

The Voronoi diagram constructed for some random (rough) approximate solution of the problem of packing $k$ equal circles  of maximum diameter in $\Omega_i$ was used as a zero estimate. The corresponding optimization problem is written as follows:
\begin{equation}
\label{opt2}
\psi(V) = \min \left\{ \min_{p, t}  \{ (v_p, c_t) + b_t \}, \; \frac{1}{2}\min_{p,q}\{ \| v_p - v_q \| \} \right\} \rightarrow \max.    
\end{equation}

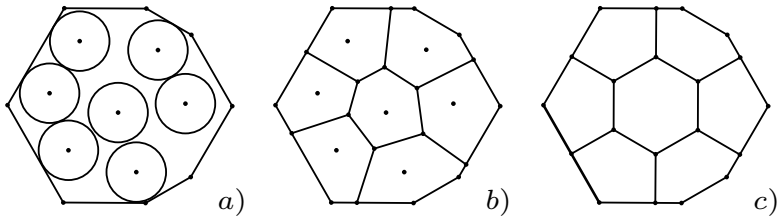
\begin{figure}
    \centering
    \begin{tabular}{p{3cm} p{3cm}  p{3cm}} 
    \definecolor{wwwwww}{rgb}{0.1,0.1,0.1}
\begin{tikzpicture}[line cap=round,line join=round,>=triangle 45,x=1.0cm,y=1.0cm, scale = 0.22]
\draw (7.6,-12.72)-- (4.23,-6.83);
\draw (4.23,-6.83)-- (7.66,-0.97);
\draw (7.6,-12.72)-- (12.57,-12.75);
\draw (12.57,-12.75)-- (15.3,-11.18);
\draw (15.3,-11.18)-- (17.81,-6.89);
\draw (17.81,-6.89)-- (15.34,-2.58);
\draw (15.34,-2.58)-- (12.62,-0.99);
\draw (7.66,-0.97)-- (12.62,-0.99);
\draw(10.91,-7.28) circle (1.8cm);
\draw(8.59,-2.96) circle (1.8cm);
\draw(6.77,-6.11) circle (1.8cm);
\draw(7.91,-9.56) circle (1.8cm);
\draw(12.01,-10.88) circle (1.8cm);
\draw(14.98,-6.74) circle (1.8cm);
\draw(13.32,-3.5) circle (1.8cm);
\begin{scriptsize}
\fill   (4.23,-6.83) circle (4pt);
\fill   (7.6,-12.72) circle (4pt);
\fill   (17.81,-6.89) circle (4pt);
\fill   (7.66,-0.97) circle (4pt);
\fill   (12.62,-0.99) circle (4pt);
\fill   (15.34,-2.58) circle (4pt);
\fill   (15.3,-11.18) circle (4pt);
\fill   (12.57,-12.75) circle (4pt);
\fill   (10.91,-7.28) circle (4pt);
\fill   (8.59,-2.96) circle (4pt);
\fill   (6.77,-6.11) circle (4pt);
\fill   (7.91,-9.56) circle (4pt);
\fill   (12.01,-10.88) circle (4pt);
\fill   (14.98,-6.74) circle (4pt);
\fill   (13.32,-3.5) circle (4pt);
\draw (16.19,-11.13) node[anchor=north west] {$ a) $};
\end{scriptsize}
\end{tikzpicture} & \begin{tikzpicture}[line cap=round,line join=round,>=triangle 45,x=1.0cm,y=1.0cm, scale = 0.22]
\draw (7.6,-12.72)-- (4.23,-6.83);
\draw (4.23,-6.83)-- (7.66,-0.97);
\draw (7.6,-12.72)-- (12.57,-12.75);
\draw (12.57,-12.75)-- (15.3,-11.18);
\draw (15.3,-11.18)-- (17.81,-6.89);
\draw (17.81,-6.89)-- (15.34,-2.58);
\draw (15.34,-2.58)-- (12.62,-0.99);
\draw (7.66,-0.97)-- (12.62,-0.99);
\draw (11.21,-0.98)-- (10.82,-4.54);
\draw (10.82,-4.54)-- (9.2,-5.41);
\draw (9.2,-5.41)-- (8.64,-7.41);
\draw (8.64,-7.41)-- (10.2,-9.46);
\draw (10.2,-9.46)-- (13.15,-8.56);
\draw (13.15,-8.56)-- (12.78,-5.8);
\draw (10.82,-4.54)-- (12.78,-5.8);
\draw (12.78,-5.8)-- (16.2,-4.08);
\draw (13.15,-8.56)-- (15.75,-10.42);
\draw (10.2,-9.46)-- (9.15,-12.73);
\draw (8.64,-7.41)-- (5.21,-8.53);
\draw (6.1,-3.62)-- (9.2,-5.41);
\begin{scriptsize}
\fill    (4.23,-6.83) circle (4pt);
\fill    (7.6,-12.72) circle (4pt);
\fill    (17.81,-6.89) circle (4pt);
\fill    (7.66,-0.97) circle (4pt);
\fill    (12.62,-0.99) circle (4pt);
\fill    (15.34,-2.58) circle (4pt);
\fill    (15.3,-11.18) circle (4pt);
\fill    (12.57,-12.75) circle (4pt);
\fill    (10.91,-7.28) circle (4pt);
\fill    (8.59,-2.96) circle (4pt);
\fill    (6.77,-6.11) circle (4pt);
\fill    (7.91,-9.56) circle (4pt);
\fill    (12.01,-10.88) circle (4pt);
\fill    (14.98,-6.74) circle (4pt);
\fill    (13.34,-3.48) circle (4pt);
\fill    (8.64,-7.41) circle (4pt);
\fill    (10.2,-9.46) circle (4pt);
\fill    (13.15,-8.56) circle (4pt);
\fill    (12.78,-5.8) circle (4pt);
\fill    (10.82,-4.54) circle (4pt);
\fill    (9.2,-5.41) circle (4pt);
\fill    (6.1,-3.62) circle (4pt);
\fill    (5.21,-8.53) circle (4pt);
\fill    (9.15,-12.73) circle (4pt);
\fill    (15.75,-10.42) circle (4pt);
\fill    (11.21,-0.98) circle (4pt);
\fill    (16.2,-4.08) circle (4pt);
\draw (16.19,-11.13) node[anchor=north west] {$ b) $};
\end{scriptsize}
\end{tikzpicture} & \definecolor{uququq}{rgb}{0.1,0.1,0.1}
\begin{tikzpicture}[line cap=round,line join=round,>=triangle 45,x=1.0cm,y=1.0cm, scale = 0.22]
\draw [line width=1pt] (7.6,-12.72)-- (4.23,-6.83);
\draw (4.23,-6.83)-- (7.66,-0.97);
\draw (7.6,-12.72)-- (12.57,-12.75);
\draw (12.57,-12.75)-- (15.3,-11.18);
\draw (15.3,-11.18)-- (17.81,-6.89);
\draw (17.81,-6.89)-- (15.34,-2.58);
\draw (15.34,-2.58)-- (12.62,-0.99);
\draw (7.66,-0.97)-- (12.62,-0.99);
\draw (11.03,-3.92)-- (8.48,-5.38);
\draw (8.48,-5.38)-- (8.47,-8.32);
\draw (8.47,-8.32)-- (11.01,-9.8);
\draw (13.56,-8.34)-- (13.57,-5.4);
\draw (13.57,-5.4)-- (11.03,-3.92);
\draw (11.01,-9.8)-- (13.56,-8.34);
\draw (11.05,-0.98)-- (11.03,-3.92);
\draw (8.48,-5.38)-- (5.94,-3.9);
\draw (8.47,-8.32)-- (5.92,-9.77);
\draw (11.01,-9.8)-- (10.99,-12.74);
\draw (13.56,-8.34)-- (16.1,-9.82);
\draw (13.57,-5.4)-- (16.12,-3.95);
\begin{scriptsize}
\fill    (4.23,-6.83) circle (4pt);
\fill    (7.6,-12.72) circle (4pt);
\fill    (17.81,-6.89) circle (4pt);
\fill    (7.66,-0.97) circle (4pt);
\fill    (12.62,-0.99) circle (4pt);
\fill    (15.34,-2.58) circle (4pt);
\fill    (15.3,-11.18) circle (4pt);
\fill    (12.57,-12.75) circle (4pt);
\fill    (5.94,-3.9) circle (4pt);
\fill    (5.92,-9.77) circle (4pt);
\fill    (11.05,-0.98) circle (4pt);
\fill    (10.99,-12.74) circle (4pt);
\fill    (16.12,-3.95) circle (4pt);
\fill    (16.1,-9.82) circle (4pt);
\fill    (11.03,-3.92) circle (4pt);
\fill    (8.48,-5.38) circle (4pt);
\fill    (8.47,-8.32) circle (4pt);
\fill    (11.01,-9.8) circle (4pt);
\fill    (13.56,-8.34) circle (4pt);
\fill    (13.57,-5.4) circle (4pt);
\draw (16.19,-11.13) node[anchor=north west] {$ c) $};
\end{scriptsize}
\end{tikzpicture}\\
    \end{tabular}
    
    \caption{A rough approximation for dense circle packing (a), the Voronoi diagram (b), and the final partition (c).}
    \label{alg_example}
\end{figure}

Here we assume that the centers of the circles are  given by the set $V= \{v_1, \dots , v_k \}$. The generation of the Voronoi diagram is performed many times, and in each case a local minimum in the problem (\ref{opt1}) is computed. The sequence of computations is shown in Fig. \ref{alg_example}.

In cases where the global minimum in the problem (\ref{opt1}) is known, the presented Algorithm \ref{alg_part} finds it rather quickly. On the other hand, in the general case we cannot claim that the optimal partition will be found with positive probability.

\begin{algorithm}[h!]

\caption{Stochastic search for the sub-optimal partition of a polygon (or a polyhedron) $\Omega$}
\label{alg_part}

\begin{algorithmic}
\STATE Input: $\Omega$, $k$

\FOR{$1 \leqslant i \leqslant N$}

\STATE $V_0 = \{v_{1}, \dots , v_k \}$ are random points distributed uniformly in $[-1,1]^n$;

\COMMENT{Run $s$ steps of Adam optimizer for the problem (\ref{opt2})}
\FOR{$1 \leqslant j < s$}
 \STATE $V_{j+1} \gets \operatorname{Adam}\left(\psi(V)\to \max ; V_j\right)$;
\ENDFOR

\COMMENT{Initialize a set of partition vertices with a Voronoi diagram}

\STATE $X_0 \gets \operatorname{Vor}(\Omega; V_s)$;
\STATE $t \gets 0$;

\COMMENT{Find the local minimum for the problem (\ref{opt1})}
\WHILE{not $\operatorname{OptCondition}(X)$}
\STATE $X_{t+1} \gets \operatorname{Adam}\left(\phi(X) \to \min; X_t\right)$;
\STATE $t \gets t+1$;
\ENDWHILE
\IF{ $\phi(X_t) < \phi(X^*)$} 
\STATE $X^* \gets \Tilde{X_t}$;
\ENDIF

\ENDFOR

\STATE Output: $X^*$.
\end{algorithmic}

\end{algorithm}

The application at the end of the article shows exactly how the partitions into $k$ parts of all sets from the UCS $S_{10}$ look like  for $10 \leqslant k \leqslant 17$. It should be noted that in the three-dimensional case the chance of obtaining the desired initial approximation is quite small. Namely, about $3\cdot 10^4$ of runs  were required to find the examples shown in Fig. \ref{fig:TRD1}, \ref{fig:TRD2}.

\section{Lower bounds}

\subsection {General scheme}

We prove a lower bound for a circle $D$ of unit diameter. When an arbitrary circle is covered by $k$ sets, there are two types of sets, namely, those that have at least two common points with the boundary of the circle (\textit{extreme sets}), and the rest (\textit{central sets}). Let there be $e$ of extreme and $c$ of central sets. Obviously, $c + e = k$.

The general scheme of the proof consists in analyzing the cases of the number of central and extreme sets, as a rule, only those where $c = 1$ or $c = 2$ remain meaningful from these cases. In each of these cases, the length of a certain segment $Q_2Q_j$ in the central set is estimated. The estimate is proved by introducing parameters-angles $\alpha, \beta, \gamma$, that is, the length of the segment $Q_2Q_j$ is considered as a function $f$ of $(\alpha, \beta, \gamma)$. The first one proves that for a fixed $\gamma$ and a fixed parameter $\delta = \gamma + 2\alpha + 2\beta$, the minimum of the function $f$ is achieved for $\alpha = \beta$. Further, it is proved that for $ \alpha = \beta$(and a fixed $ \delta$) , the greater the $\gamma$ , the smaller the value of $f$, so for the lower estimate, we need to take the maximum possible value of $\gamma = \gamma_{\max}$. Finally, it is proved that for $\gamma = \gamma_{\max}, \alpha = \beta$, the function $f$ is minimal for the minimum possible value of $\delta$, that is, for the minimum possible values of $\alpha = \beta = \alpha_{\min}$. As a result, the length of the segment $Q_2Q_j$ can be estimated from below by the value $f(\alpha, \beta, \gamma) = f (\alpha_{\min}, \alpha_{\min}, \gamma_{\max})$.

Below is a table of parameters for various cases of $k$, $c$.

\setlength{\tabcolsep}{9pt}
\renewcommand{\arraystretch}{0.95}
\begin{center}
\begin{table}[H]
\caption{Table of extreme parameter values}
\label{results2}
 \begin{tabular}{ |c|c|c|c|c|c|c|c|  }
 \hline
 $k$ & $c$ & $e$ & $j(k, c)$ & $\alpha_{\min}$ & $\gamma_{\max}$ & $\delta_0$ & 
$Q_2 Q_j$ \\ 
 \hline\hline
 10 & 1 & 9  & 6 &  $36^{\circ}$ & $86.92^{\circ}$      &  $230.92^{\circ}$    & $0.3665$ \\
 \hline
 11 & 1 & 10 & 7 &  $27.88^{\circ}$ & $124.23^{\circ}$ & $235.79^{\circ}$ & $0.3535$ \\ 
 \hline
 11 & 2 & 9 &  6 & $38.23^{\circ}$ & $82.82^{\circ}$ & $235.79^{\circ}$ & $0.3665$\\
 \hline
 12 & 1 & 11 & 7 & $20^{\circ}$ & $120^{\circ}$ & $200^{\circ}$ & $\sin(\frac{\pi}{9})$ \\
 \hline
 12 & 2 & 10 & 6 & $30^{\circ}$ & $80^{\circ}$ & $200^{\circ}$ & $0.3751$ \\
 \hline

\end{tabular}
\end{table}
\end{center}

\subsection {Complete proof}

Define  $\sigma_{10} = 0.366538$, $\sigma_{11} = 0.353553$, $\sigma_{12} = \sin(\frac{\pi}{9}) = 0.342020...$ We prove that $d_{k}(D) \geqslant \sigma_k$ for $k \in \{10, 11, 12\}$.

For $c = 0$, the set covering the center of the circle will have a diameter of at least $0.5$, so we do not consider this case. For $c \geqslant 3$, we have $e \leqslant k-3$, which means that a diameter of at least  $\sin(\frac{\pi}{k - 3}) = \begin {cases} 0.4338... (k = 10) \\  0.3826... (k = 11) \\ 0.3420... (k = 12) \end{cases} \geqslant \sigma_k$ will be required to cover the boundary of the circle. The remaining case is $c = 1$ and $c = 2$.

Let's fix the orientation of the circle counterclockwise. For each ``extreme'' set, due to its closure, there is its first point on the circle, in accordance with the orientation. At the same time, all these first points are different. Let's denote these points (in accordance with their order when traversing the circle) by $P_1, P_2, ..., P_{e}$, and denote the corresponding sets by $S_1, S_2, ..., S_{e}$.

 Also, for the convenience of notation, we put $P_0 = P_{e}$, $P_{e + 1} = P_1$. We denote by $T_{1}, ..., T_{c}$ ``central'' sets. $T = T_{1} \cup ... \cup T_{c}$ (union of  ``central'' sets)

\begin{figure}[htb]
\center{\includegraphics[scale=0.20]{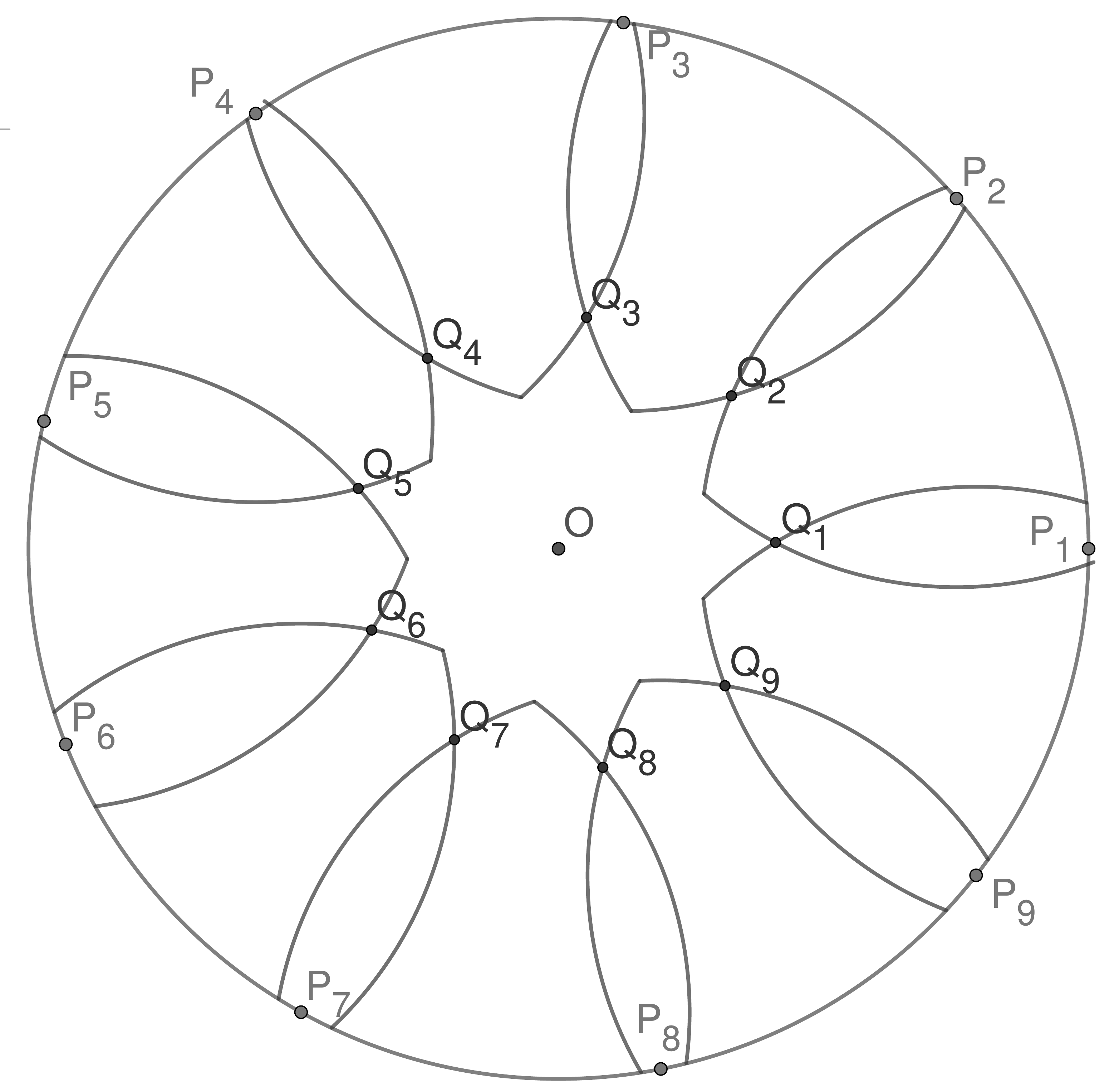}}
\caption{To the proof of the lower bound $d_{10}$, $d_{11}, d_{12}$ (example for $e = 9$)}
\label{lowerA}
\end{figure}
 
Denote by $Q_i$ the common point of $D$ and circles with centers $P_{i-1}$, $P_{i + 1}$ of radius $\sigma_k$ (see Figure \ref{lowerA}). Due to the closeness of the sets, we can assume that $Q_i \in T$.

 \vspace{\baselineskip}

Our goal is to prove the following fact $$\operatorname{dist}(Q_2, Q_{j}) \geqslant \sigma_k \label{fact1}$$

where $j = j(k, c) = \begin{cases}
6, \ \ k=10, c=1\\
7, \ \ k=11, c=1\\
6, \ \ k=11, c=2 \\
7, \ \ k=12, c=1 \\
6, \ \ k=12, c=2
\end{cases}$

 \vspace{\baselineskip}

\begin{lemma}
 $Q_2$ and $Q_j$ lie in the same central set up to renumbering.
 \label{centralset}
\end{lemma}

\textbf{Proof:} In the case of $c = 1$, all points $Q_1,..., Q_e$ belong to the central set, so the lemma is valid.
In the case of $c = 2$, some points belong to $T_1$, some to $T_2$. We show that it is possible to renumber the points so that $Q_2$ and $Q_j$ belong to the same central set.

For $k = 11$, $c = 2$, we have $j = 6$. Note that the graph $G = \{V = (1,..., 9), E = \{i, j\, |\, i - j = 4 \pmod 9\}$ is not bipartite, which means that for any division into two sets there will be two numbers $x, y$ such that $x - y = 4 \pmod 9$. Renumber the vertices so that $y = 2$, $x = 6 = j$.

For $k = 12$, $c = 2$, we have $j = 6$. The graph $G = \{V = (1,..., 10), E = \{i, j \,|\, i - j = 4 \pmod {10} \}$ is not bipartite, which means that for any division into two sets there will be two numbers $x, y$ such that $x - y = 4 \pmod 9$. Renumber the vertices so that $y = 2$, $x = 6 = j$. \qed

 \vspace{\baselineskip}

It follows from Lemma \ref{centralset} and (\ref{fact1}) that \[ d_{k} (D) \geqslant \max(\operatorname{diam}(T_1), ..., \operatorname{\operatorname{diam}}(T_c)) \geqslant \operatorname{dist} (Q_2, Q_j) \geqslant \sigma_k,\] and this is what we need to prove.

 \vspace{\baselineskip}
 
The position of $ Q_2 $ is uniquely determined by the points $ P_1 $ and $P_3$, and the position of $ Q_j $ is uniquely determined by the points $ P_{j-1}$, $ P_{j + 1}$. Let $\alpha := \angle{Q_2 O P_3}$, $\beta := \angle{Q_{j} O P_{j + 1}}$, $\gamma := \angle{P_3 O P_{j-1}}$. It is obvious that $\angle{P_1 O Q_2} = \alpha$ and $\angle{Q_{j} O P_{j+1}} = \beta$. Let $\delta = \angle{P_1 O P_{j + 1}} = 2\alpha + \gamma + 2\beta$.

\begin{figure}[htb]
\center{\includegraphics[scale=0.20]{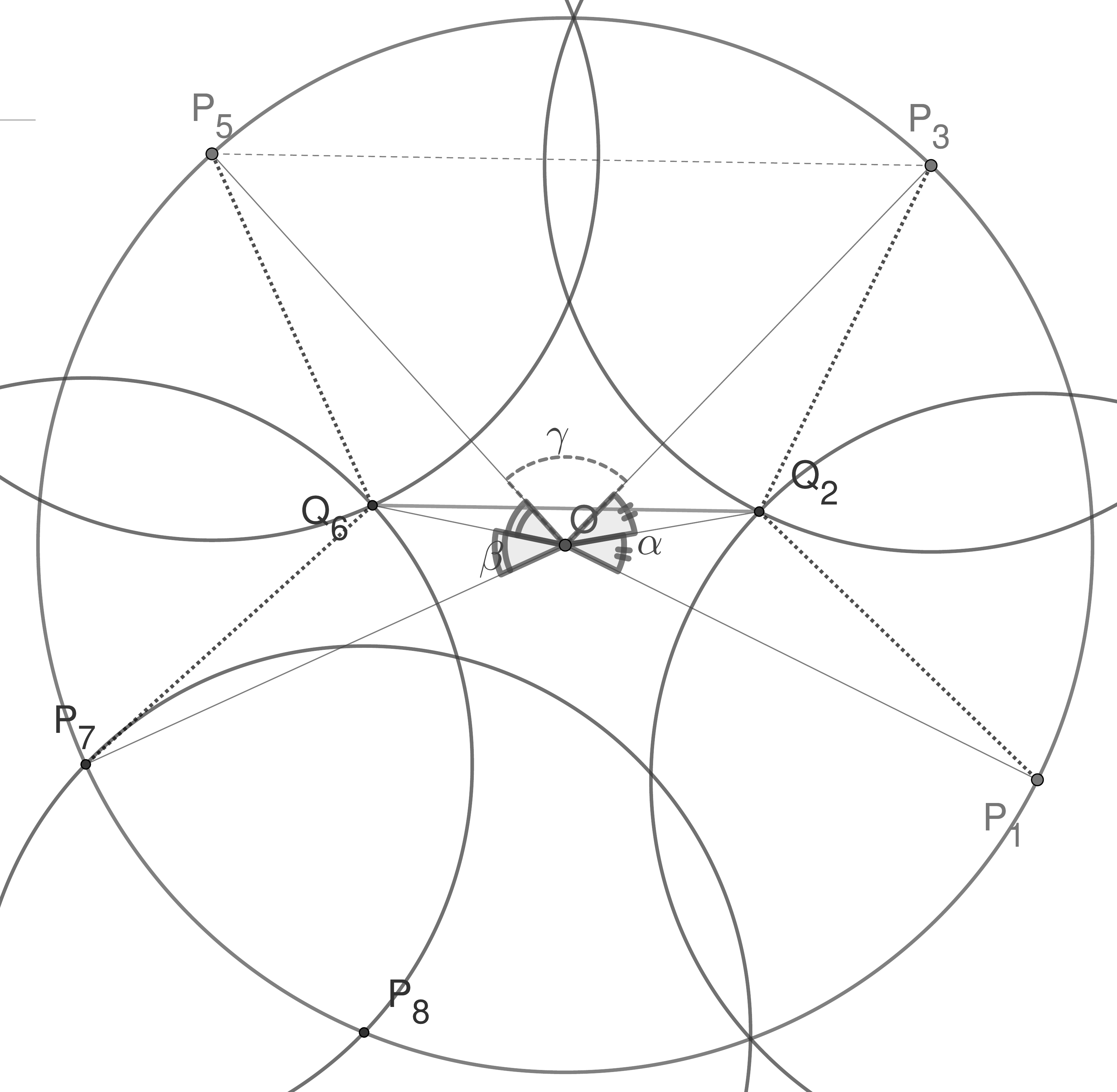}}
\caption{To the proof of the lower bound $d_{10}$, $d_{11}, d_{12}$}
\label{lowerB}
\end{figure}

Note that $\operatorname{dist}(Q_2, Q_{j})$ is a function of $f (\alpha, \beta, \gamma)$ that depends only on the angles $\alpha$, $\beta$, $\gamma$. We show that this function reaches a minimum when $\alpha_0 = \beta_0 = \frac{\delta_0 - \gamma_{\max}}{4}$ and $\gamma_{\max} = (2j - 8) \cdot \operatorname{arcsin}(\rho)$.

Note that the arcs $P_{i}P_{i + 1}$ cannot be greater than $ 2 \cdot \operatorname{arcsin}(\sigma_k)$, so we have such inequalities on the angles:

$$ \gamma \leqslant \gamma_{\max} = \gamma_{\max}(k, c) =  (2j - 8) \cdot \operatorname{arcsin}(\sigma_k) = \begin{cases}
179.09^{\circ}, k=10, c=1\\
124.23^{\circ}, k=11, c=1\\
82.82^{\circ}, k=11, c=2\\
120^{\circ}, k=12, c=1 \\
80^{\circ}, k=12, c=2
\end{cases}$$ 

$$\delta \geqslant \delta_0 = \delta_0(k, c) = 2\pi - 2(k - c - j) \cdot \operatorname{arcsin}(\sigma_k) = \begin{cases}
230.92^{\circ}, k = 10,c=1\\
235.79^{\circ}, k = 11, c=1\\
235.79^{\circ}, k = 11, c=2\\
200^{\circ}, k=12, c=1\\
200^{\circ}, k=12, c=2
\end{cases}$$

 I. Let $\gamma$ and $\delta$ be fixed. We intend to show that the value of $\operatorname{dist} (Q_2, Q_j)$ is the smallest when $\alpha = \beta$. Since $\delta$ is fixed, we can assume that the positions of the points $P_1$ and $P_{j + 1}$ are fixed. Without generality restriction, let $\alpha \leqslant \beta$. We denote by $P'_3$, $P'_{j - 1}$, $Q'_2$, $Q'_{j}$ the positions of the points $P_3$, $P_{j-1}$, $Q_2$, $Q_{j}$ for $\alpha = \beta$, respectively.

\begin{figure}[htb]
\center{\includegraphics[scale=0.20]{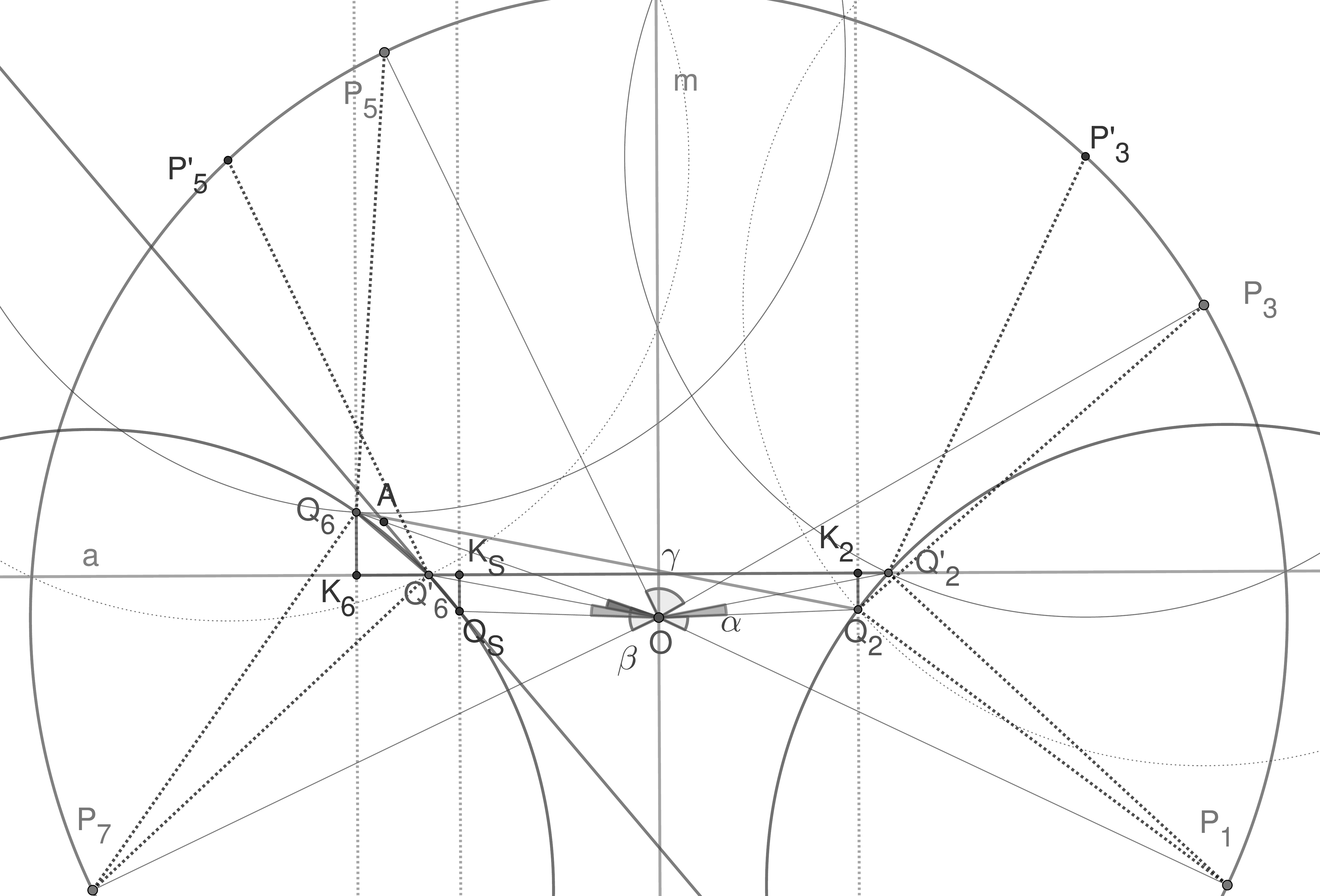}}
\caption{To the proof of the lower bound $d_{10}$, $d_{11}, d_{12}$}
\label{fig:image_lbound}
\end{figure}
 
 Denote by $a$ the line passing through $Q'_2 Q'_{j}$. Let the line $m \perp a$, $O \in m$. 
 Let $Q_S$ be a point symmetric to $Q_2$ with respect to $m$. Denote the projections of the points $Q_2$, $Q_{j}$, $Q_S$ on the line $a$ through $K_2$, $K_{j}$, $K_S$, respectively.
 
 Note that since $\gamma$ and $\delta$ are fixed, then $\angle{Q_2 O Q_{j}} = \alpha + \gamma + \beta = \frac{\gamma + \delta}{2} = \angle{Q'_2 O Q'_{j}}$. Therefore, $\angle{Q_2 O Q'_2} = \angle{Q_{j} O Q'_{j}} = \angle{Q'_{j} O Q_S}$(the latter equality is true due to the symmetry of the angles with respect to $m$).
 
 Denote $A =  O Q_{j} \cap Q_S Q'_{j} \text{(line)}$. Consider $\triangle{O Q_S A}$.In it, $OQ'_{j}$ is a bisector, and due to the limitation of $\delta$ from below in $200^{\circ}$, it can be shown that $\angle{O Q_S A} > 90^{\circ}$, which means that $Q_{j} Q'_{j} > A Q'_{j} > Q'_{j}Q_S$. We also have $\angle{K_{j} Q'_{j} Q_{j}} < \angle{Q_S Q'_{j} K_S}$. Therefore, $Q'_{j} K_S = Q'_{j}Q_S \cos(\angle{Q_S Q'_{j} K_S}) < Q'_{j} K_{j} \cos(\angle{K_{j} Q'_{j} Q_{j}}) = K_{j} Q'_{j}$. Hence, $K_2 Q'_2 = K_S Q'_{j} < Q'_{j} K_{j} $.
 
 As a result, we get the required inequality $Q_2 Q_{j} > K_2 K_{j} = K_2 Q'_{j} + Q'_{j}K_{j} > K_2Q'_{j} + K_2Q'_2 = Q'_{j} Q'_2$. This means that the value of $\operatorname{dist} (Q_2, Q_{j})$ is minimal for $\alpha = \beta$.

 \vspace{\baselineskip}
 
 II. Now let $\alpha = \beta$ and $\delta$ be fixed. Note that when $\gamma$ decreases, the value of $\angle{Q_2 O Q_{j}} = \alpha + \gamma + \beta = \frac{\gamma + \delta}{2}$ also decreases, and hence the value of $\operatorname{dist}(Q_2, Q_{j})$ increases due to the fact that the points $Q_2$ and $Q_{j}$ ``shift'' along the corresponding circles in different directions, approaching the boundary of the circle $D$ (we assume that $P_1$ and $P_{j+1}$ are fixed at this moment).
 
 On the other hand, $\gamma \leqslant \gamma_{\max} = \gamma_{\max}(k, c)$.

 \vspace{\baselineskip}
 
 III. Now let's say $\alpha = \beta$ and $\gamma = \gamma_0$. Consider the isosceles trapezoid $Q_2P_3P_{j-1}Q_{j}$. Note that when the angle $\delta$ decreases, the angles $\angle{Q_{j}P_{j-1}P_3} = \angle{P_{j-1}P_3Q_2}$ also decrease. Therefore, the length of $Q_2Q_{j}$ also decreases.
 
 On the other hand, $\delta \geqslant \delta_0(k, c)$.
 
 \vspace{\baselineskip}
 
 Summing up, we can conclude that the value of $\operatorname{dist}(Q_2, Q_{j})$ reaches its minimum at the above values of the angles. The calculations show that for the specified values $\alpha, \beta, \gamma, \delta$ we have $$f_{\min} = f(\alpha_{\min}, \alpha_{\min}, \gamma_{\max}) = \begin {cases}
 0.3665..., k=10, c=1\\
 0.3536..., k=11, c=1\\
 0.4362..., k=11, c=2\\
 \sin(\frac{\pi}{9}), k=12, c=1\\
 0.3751..., k=12, c=2\\
 \end{cases}\geqslant \sigma_k.$$ 
 Finally, some of the ``central'' sets has a diameter of at least $\operatorname{dist}(Q_2, Q_{j}) \geqslant \sigma_k$, which was required to be proved.
 
 \section*{Acknowledgements}
The authors would like to thank the anonymous reviewers for careful reading and for comments that helped improve the text of the article and correct a number of inaccuracies.
 
\bibliography{tolmachev}

\newpage
\section{Application}

\begin{figure}[htb]
\center{\includegraphics[scale=0.25]{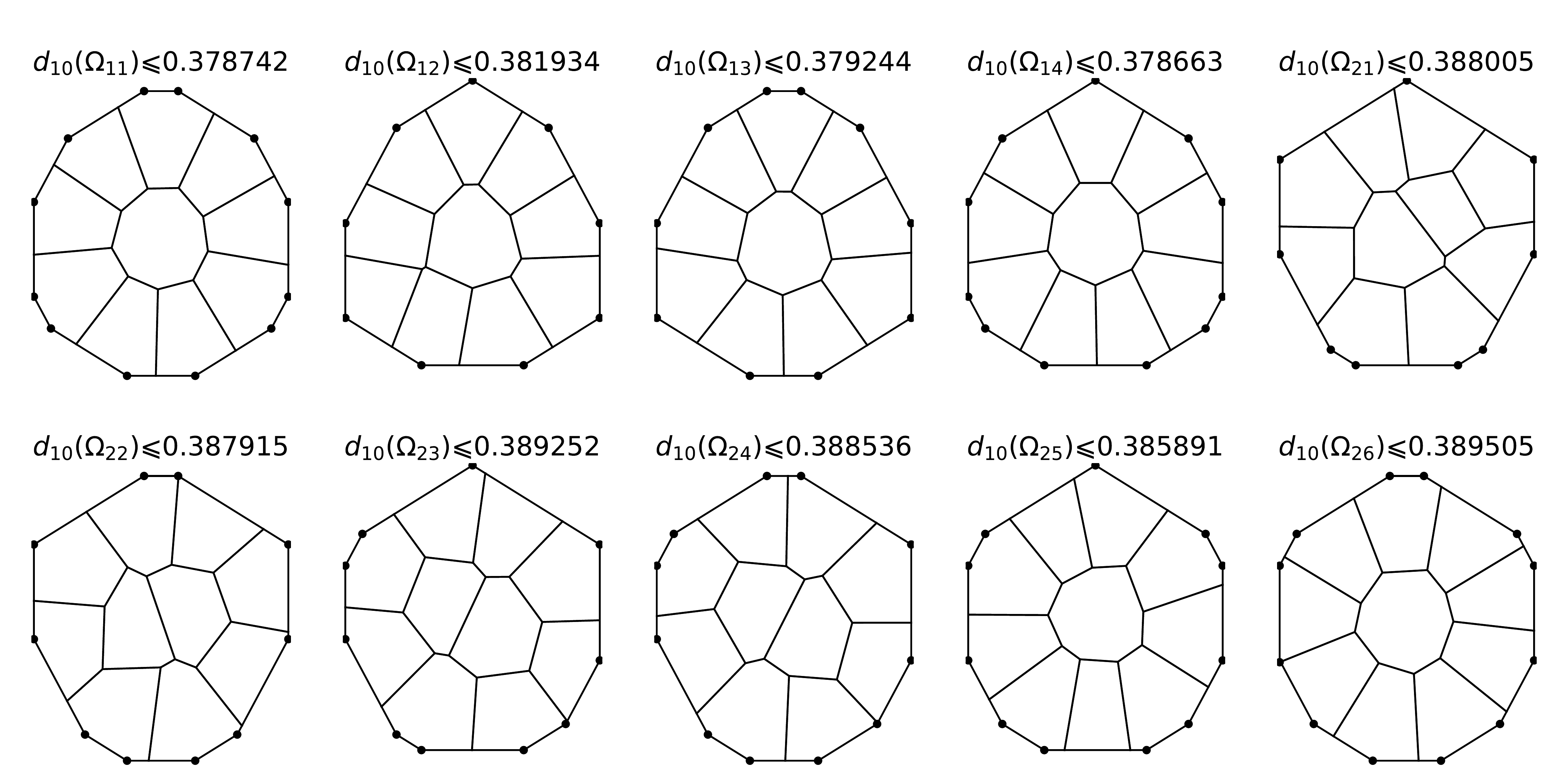}}
\caption{Splitting shapes from the UCS $S_{10}$ to improve the upper bound for $d_{10}$}
\label{fig:image1}
\end{figure}

\begin{figure}[htb]
\center{\includegraphics[scale=0.25]{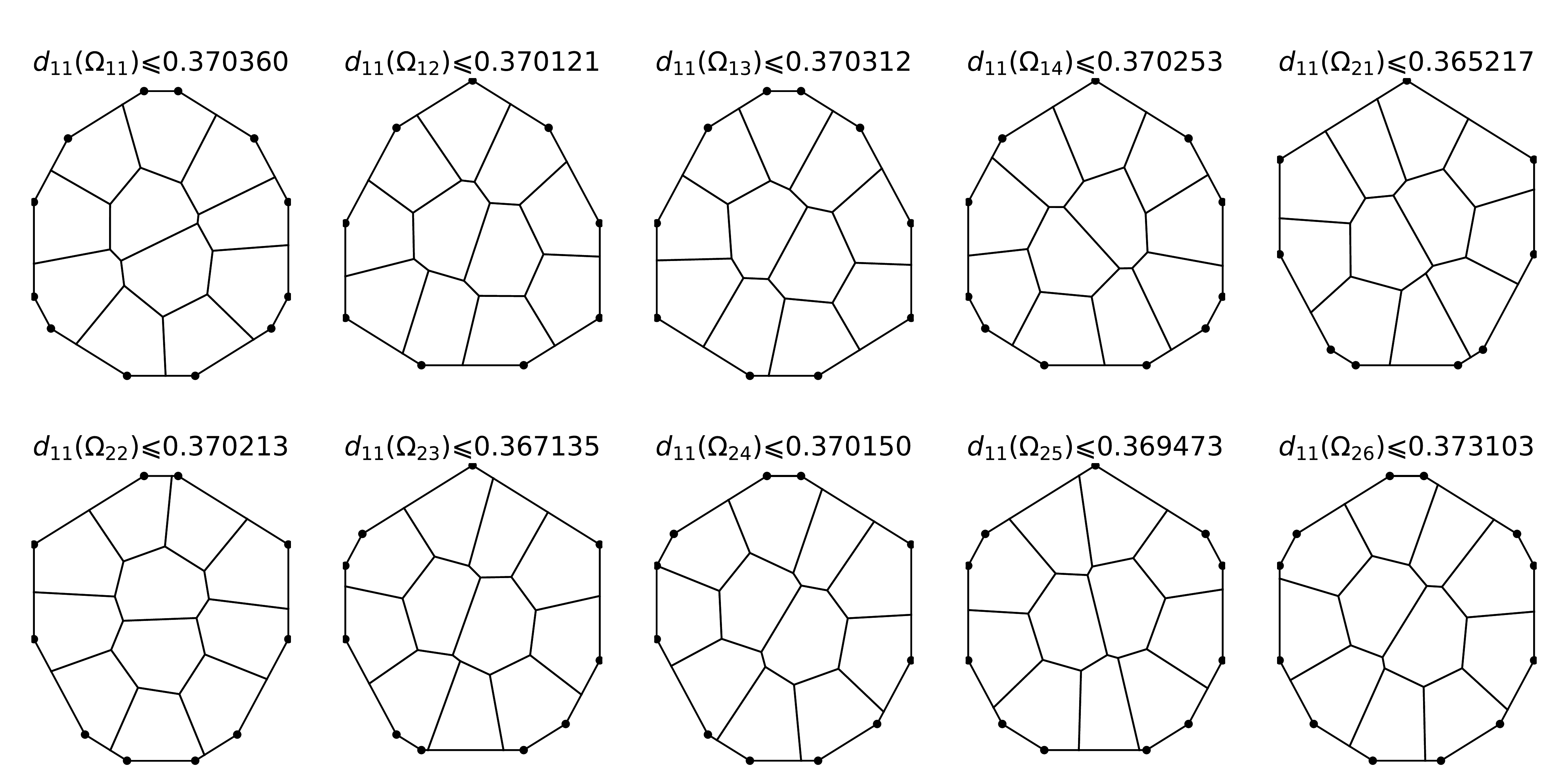}}
\caption{Splitting shapes from the UCS $S_{10}$ to improve the upper bound for $d_{11}$}
\label{fig:image2}
\end{figure}
 
\begin{figure}[htb]
\center{\includegraphics[scale=0.25]{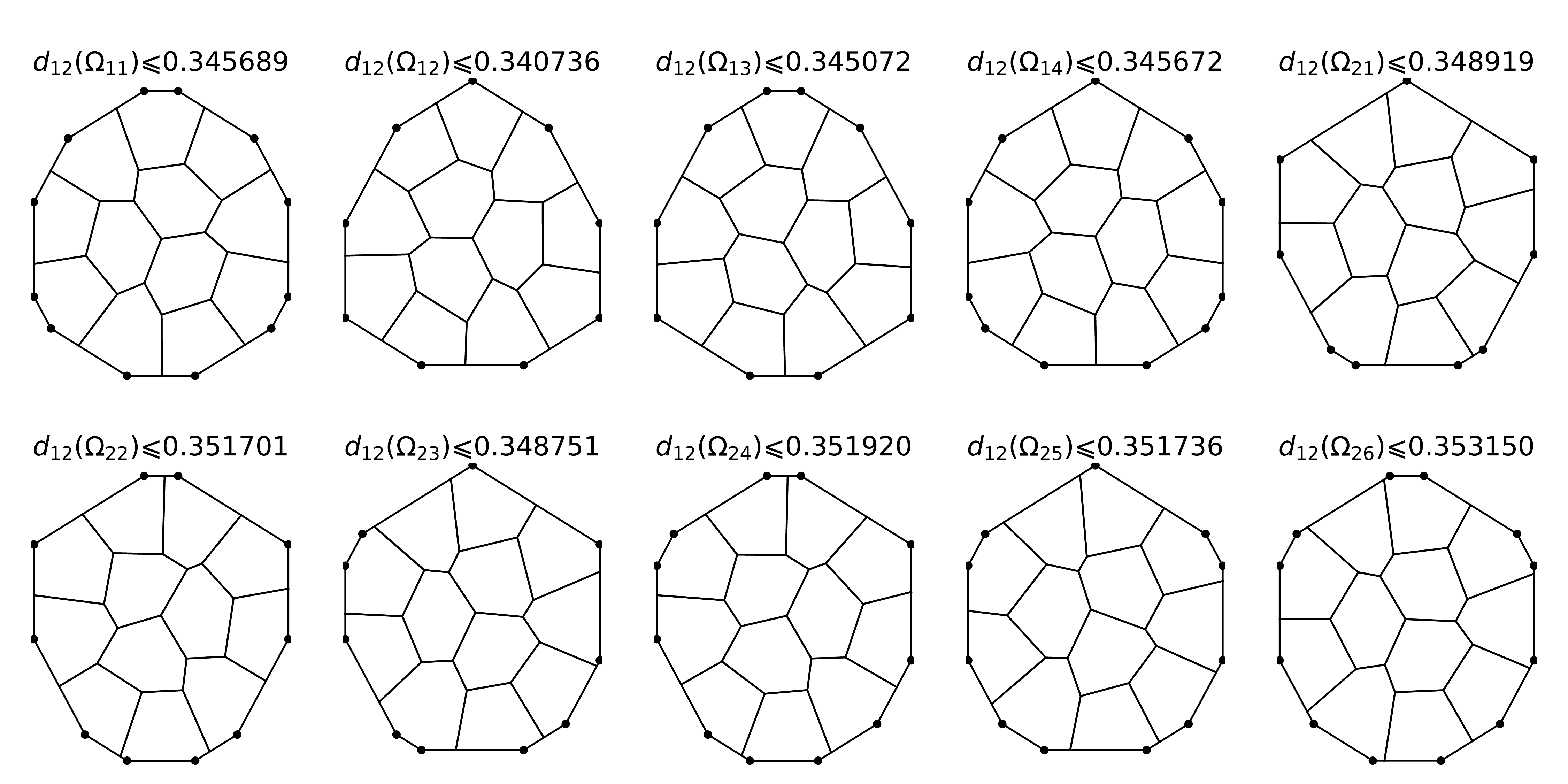}}
\caption{Splitting shapes from the UCS $S_{10}$ to improve the upper bound for $d_{12}$}
\label{fig:image3}
\end{figure}

\begin{figure}[htb]
\center{\includegraphics[scale=0.25]{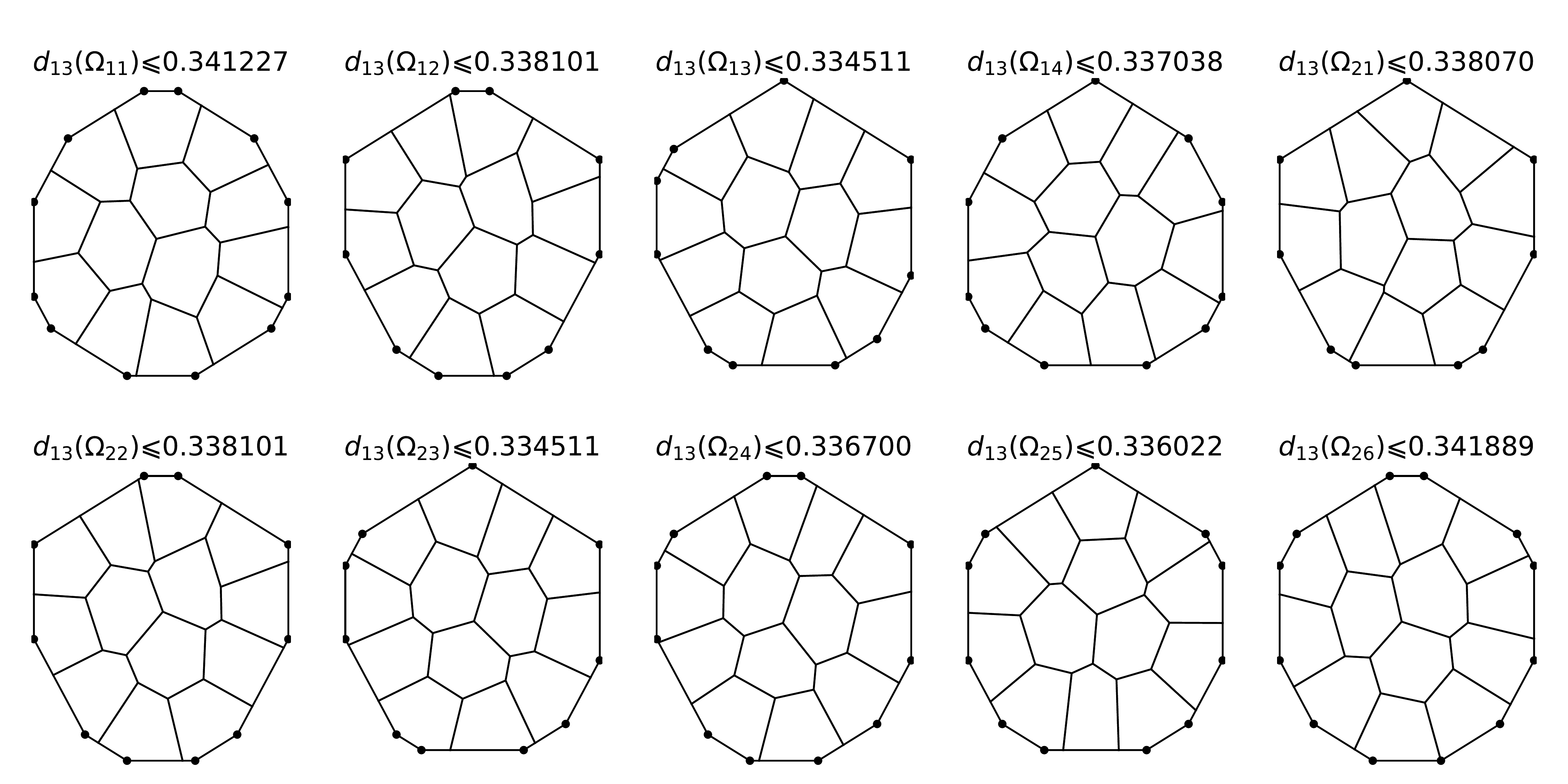}}
\caption{Splitting shapes from the UCS $S_{10}$ to improve the upper bound for $d_{13}$}
\label{fig:image4}
\end{figure}

\begin{figure}[htb]
\center{\includegraphics[scale=0.25]{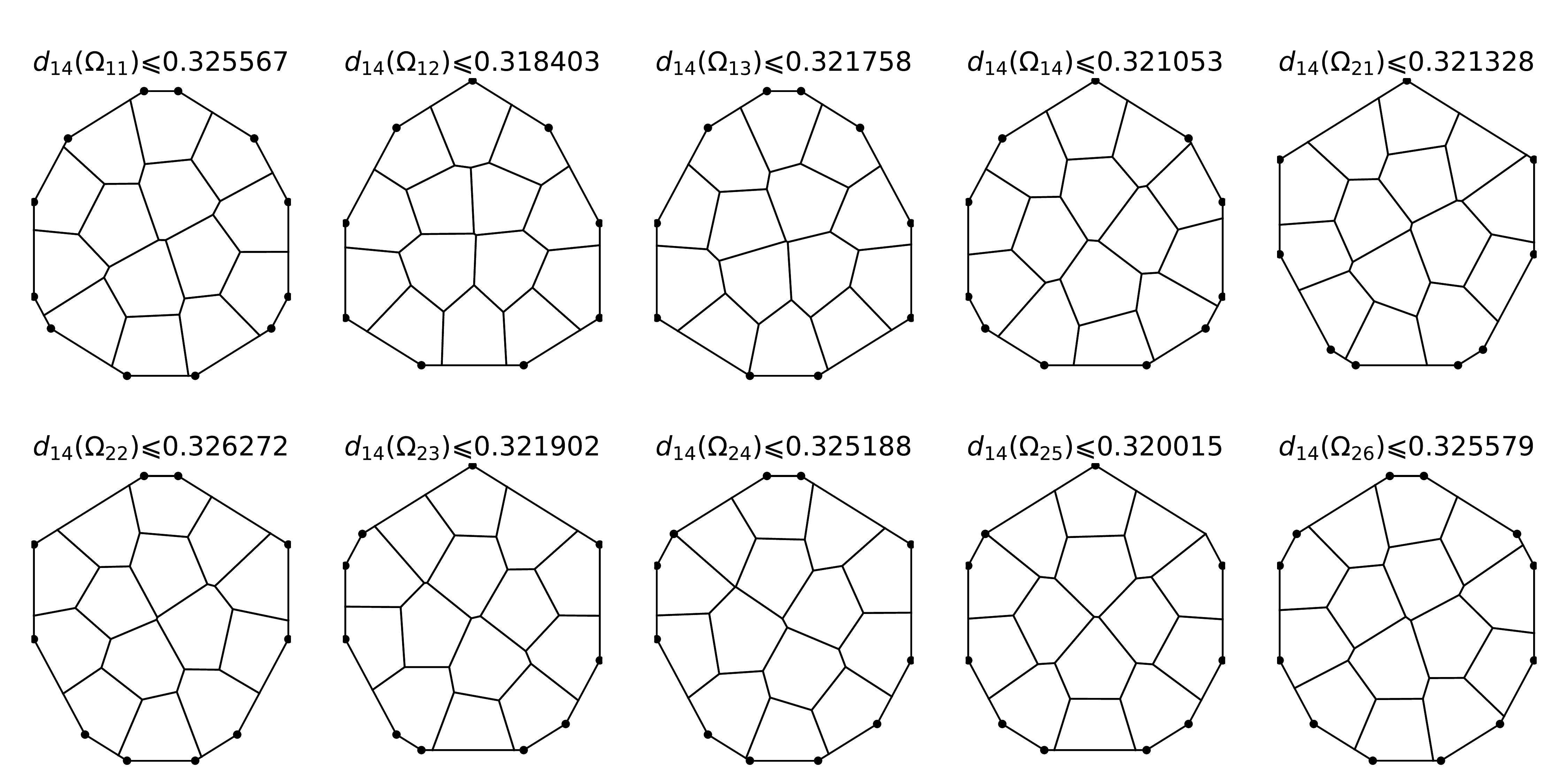}}
\caption{Splitting shapes from the UCS $S_{10}$ to improve the upper bound for $d_{14}$}
\label{fig:image5}
\end{figure}

\begin{figure}[htb]
\center{\includegraphics[scale=0.25]{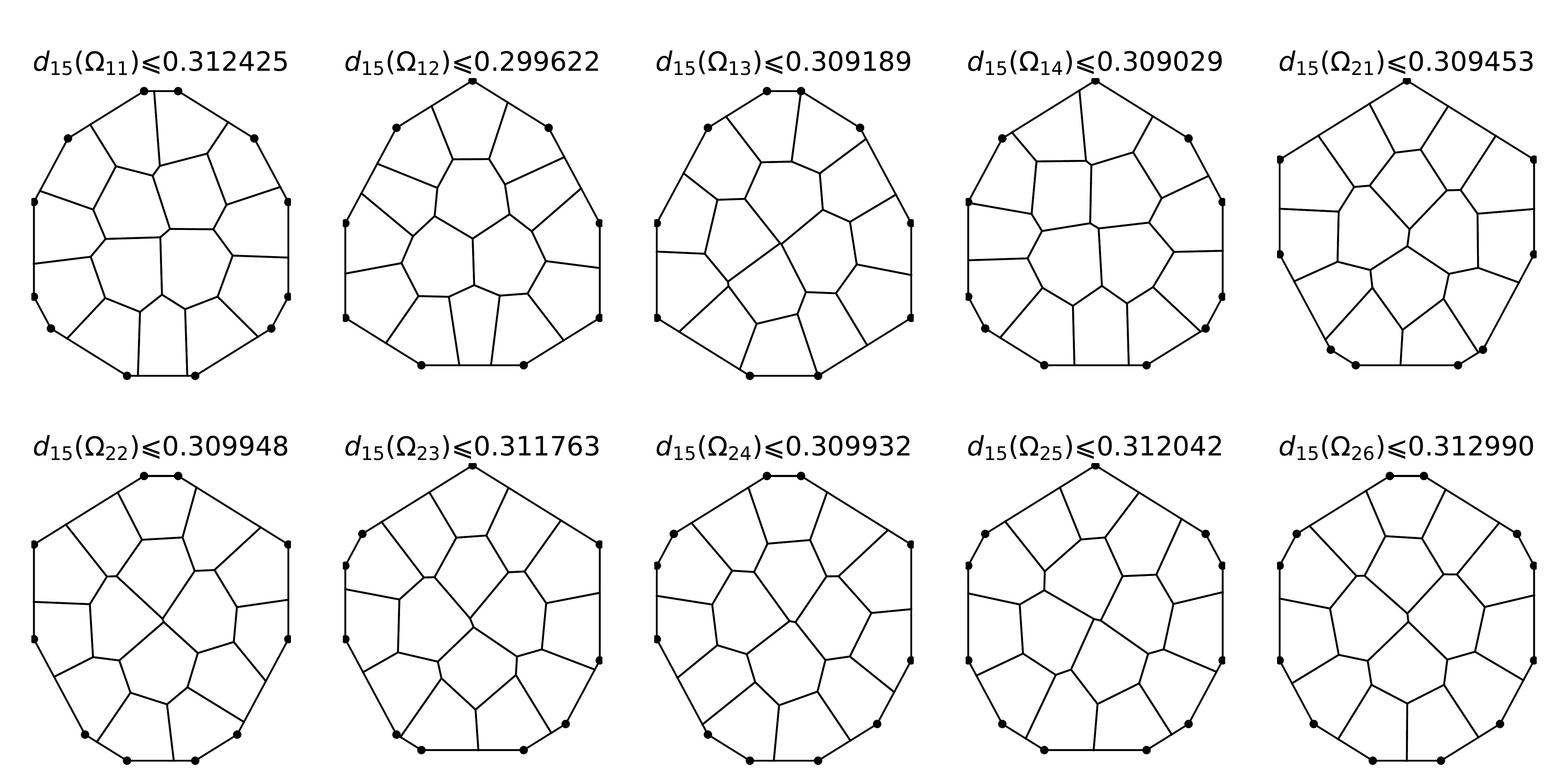}}
\caption{Splitting shapes from the UCS $S_{10}$ to improve the upper bound for $d_{15}$}
\label{fig:image6}
\end{figure}

\begin{figure}[htb]
\center{\includegraphics[scale=0.25]{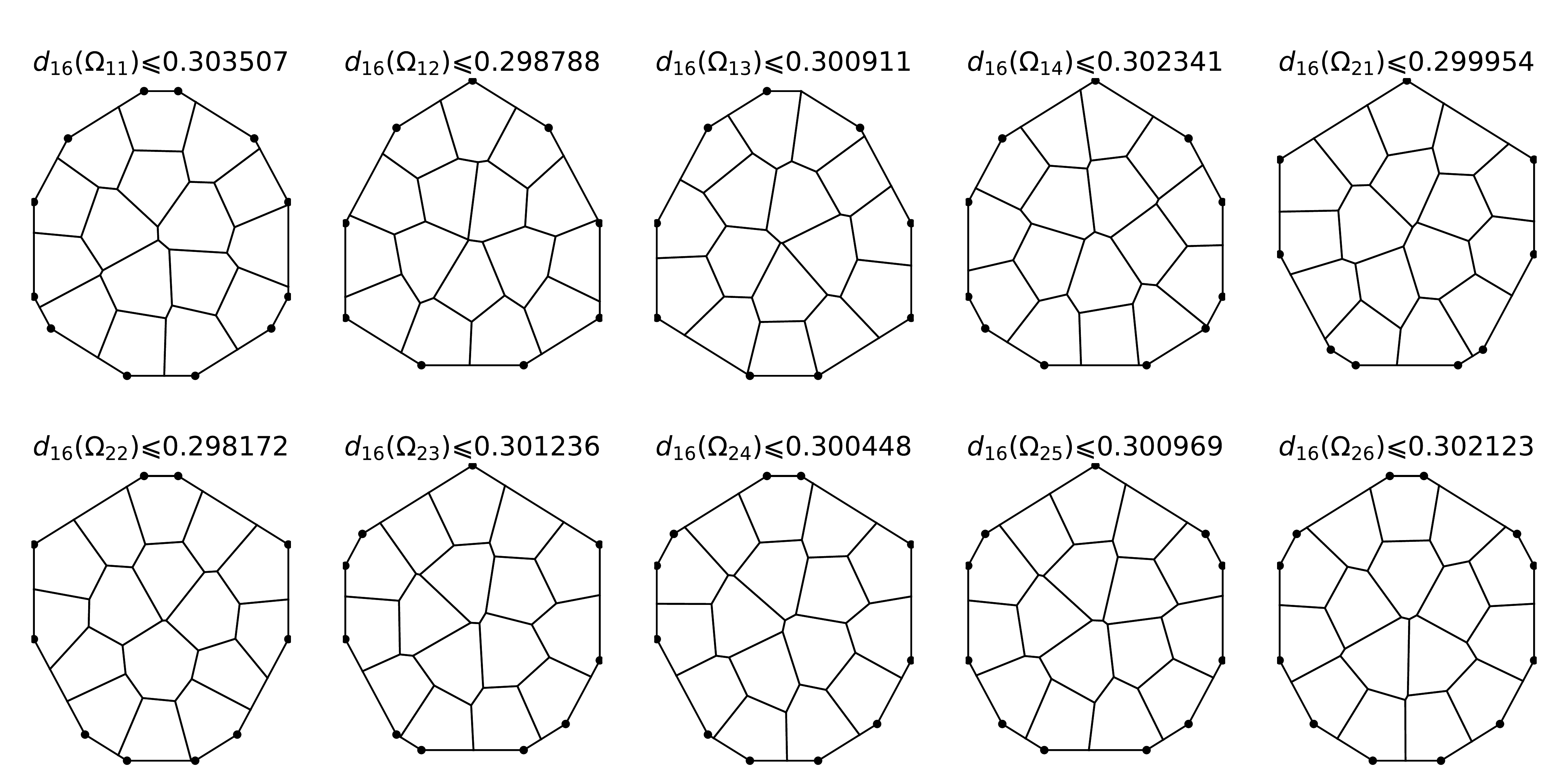}}
\caption{Splitting shapes from the UCS $S_{10}$ to improve the upper bound for $d_{16}$}
\label{fig:image7}
\end{figure}

\begin{figure}[htb]
\center{\includegraphics[scale=0.25]{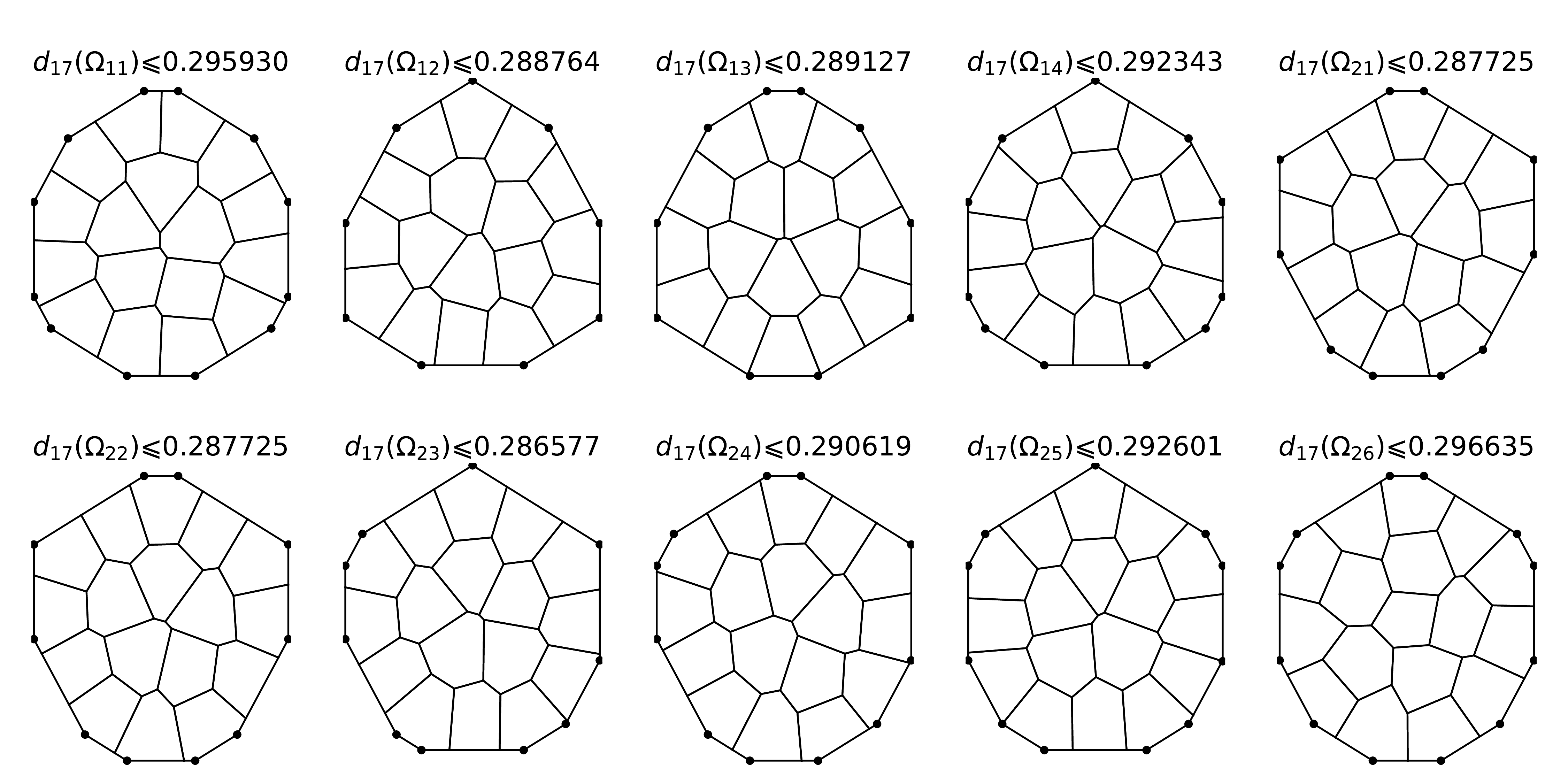}}
\caption{Splitting shapes from the UCS $S_{10}$ to improve the upper bound for $d_{17}$}
\label{fig:image8}
\end{figure}

\begin{figure}[htb]
\center{\includegraphics[scale=0.30]{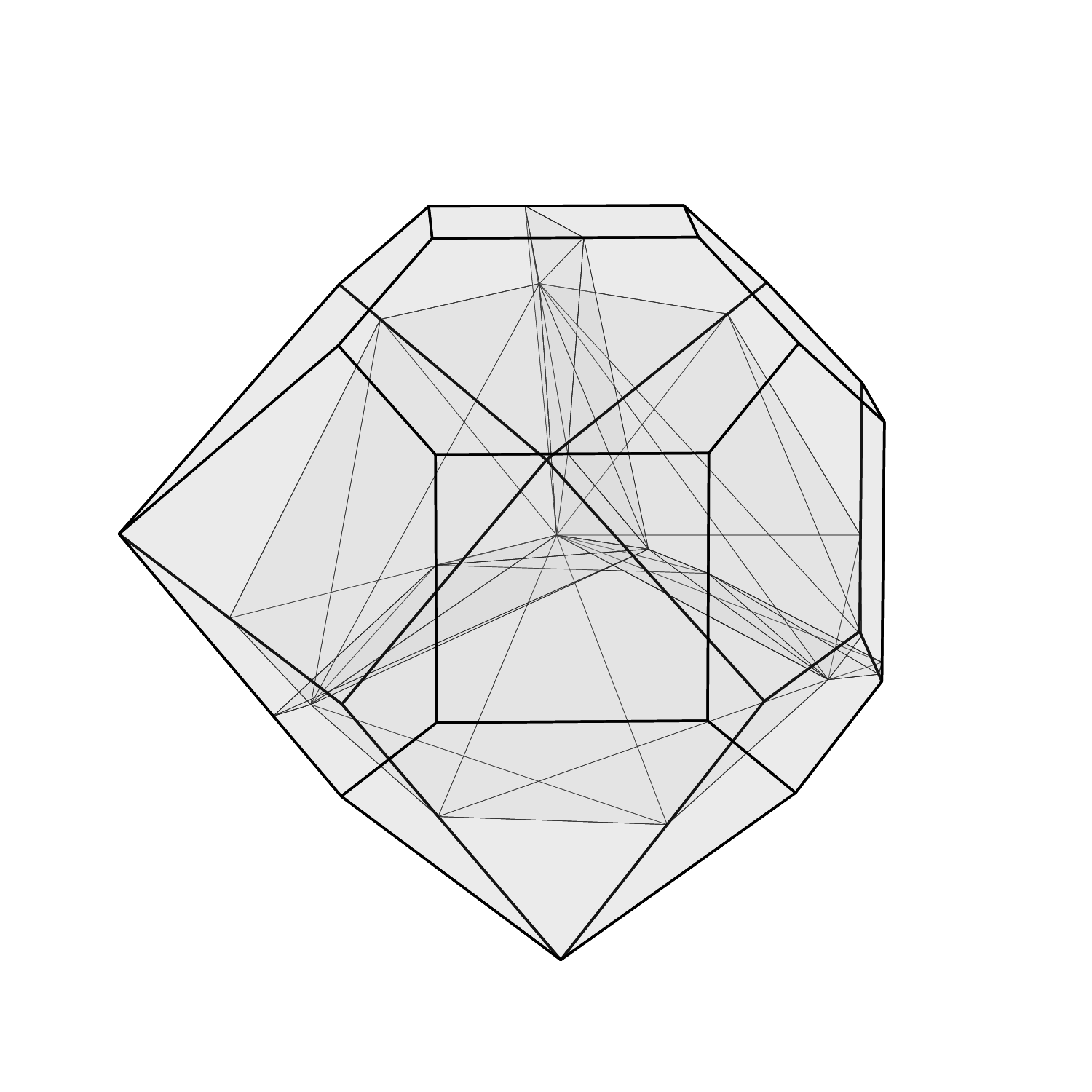} \includegraphics[scale=0.30]{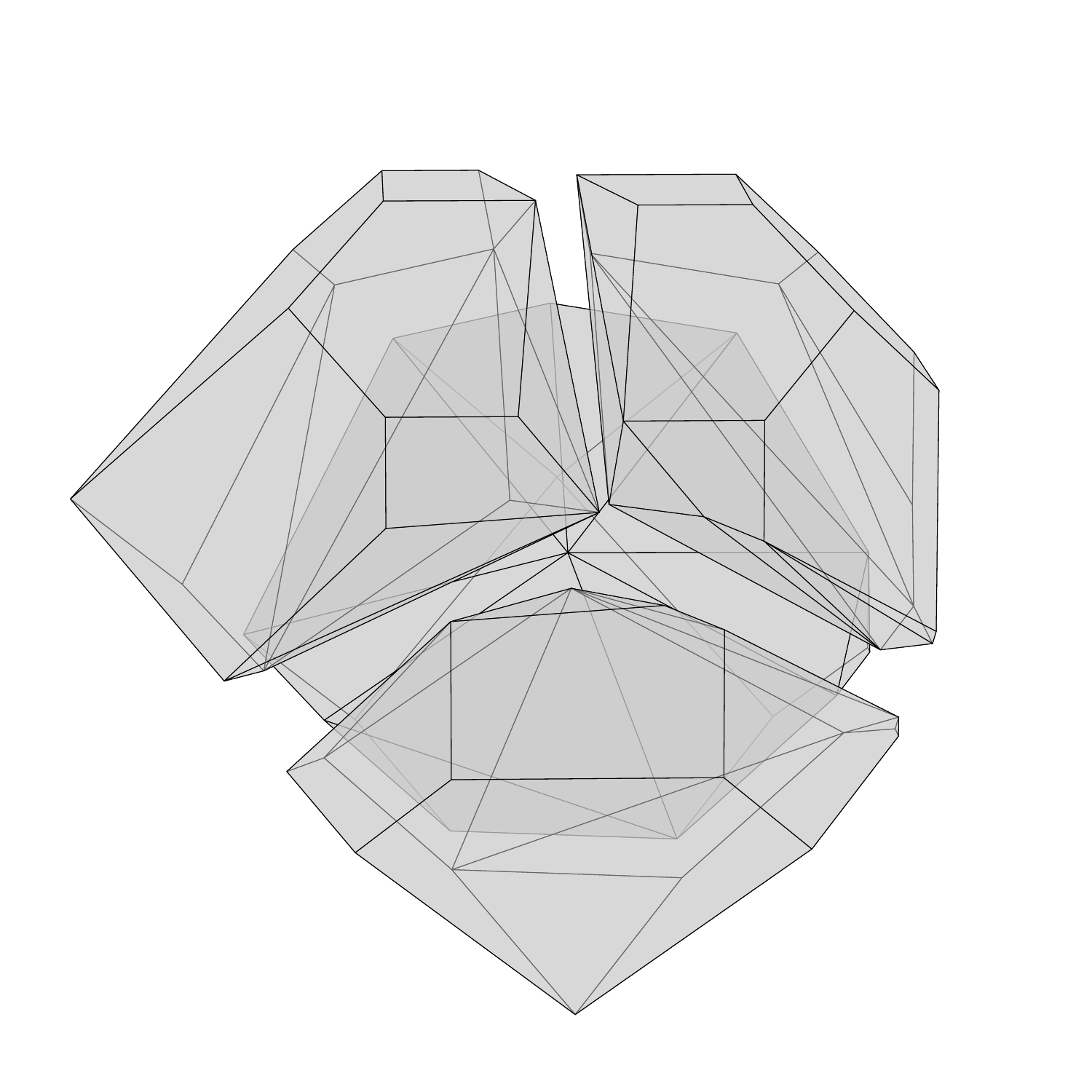}}
\caption{Covering of the truncated rhombic dodecahedron, $d_{3,4} \leqslant 0.966$}
\label{fig:TRD1}
\end{figure}

\begin{figure}[htb]
\center{\includegraphics[scale=0.30]{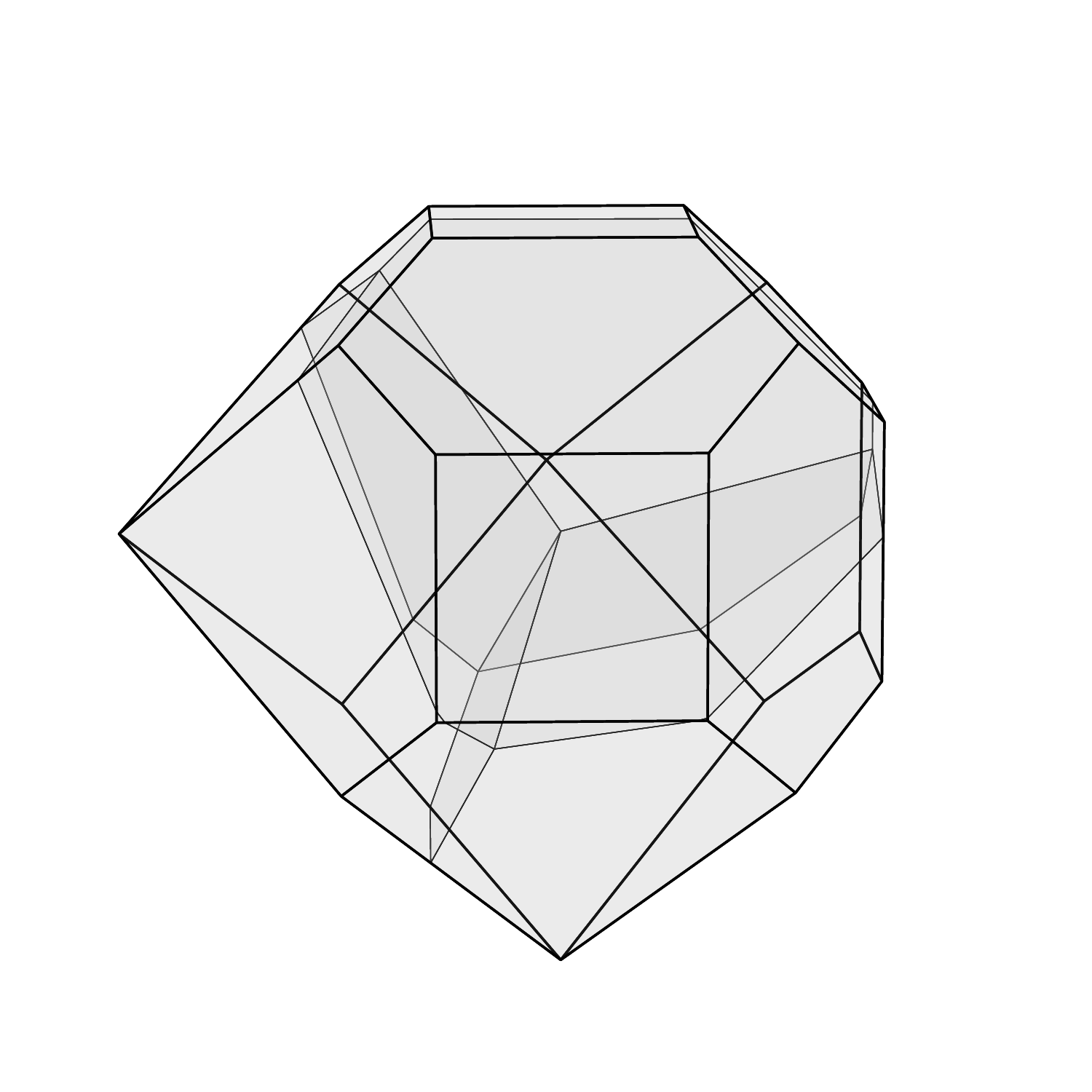} \includegraphics[scale=0.30]{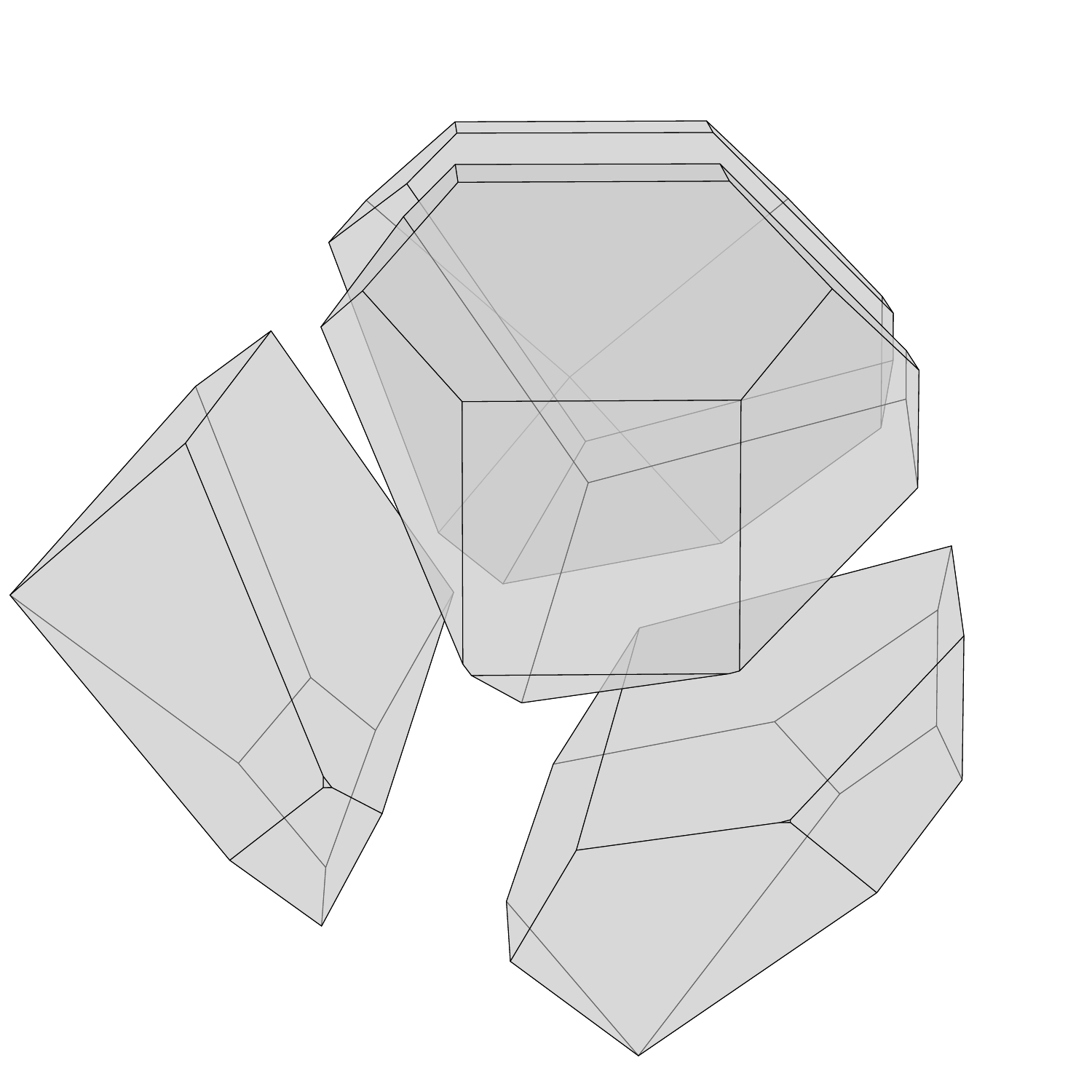}}
\caption{Partitioning  the truncated rhombic dodecahedron by six planes, $d_{3,4} \leqslant 0.9755$}
\label{fig:TRD2}
\end{figure}

\end{document}